\documentclass[reqno]{alea2}
\usepackage{natbib}
\usepackage{fancyhdr}
\usepackage{graphicx}
\eheader{Alea}{{\bf 1}}{2006}{149}{180}

\elogo{\parbox[c]{3cm}{\includegraphics[width=3cm]{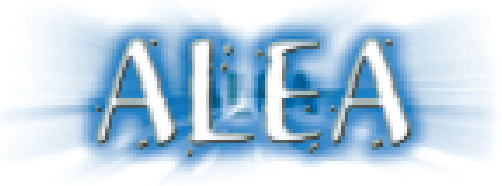}}}

\usepackage{amssymb}
\usepackage{amsmath}
\renewcommand{\cite}{\citet}

\makeatletter \@addtoreset{equation}{section} \makeatother

\renewcommand\theequation{\thesection.\arabic{equation}}
\renewcommand\thetable{\thesection.\@arabic\c@table}

\newtheorem{theorem}{Theorem}[section]
\newtheorem{lemma}[theorem]{Lemma}
\newtheorem{proposition}[theorem]{Proposition}

\newtheorem{remark}[theorem]{Remark}

\def\cL{\mathcal{L}}

\usepackage{graphics,epsf,psfrag}

\DeclareMathOperator{\sign}{\mathrm{sign}}

\newcommand{\ind}{\mathbf{1}}

\newcommand{\R}{\mathbb{R}}
\newcommand{\Z}{\mathbb{Z}}
\newcommand{\N}{\mathbb{N}}
\renewcommand{\tilde}{\widetilde}

\newcommand{\cH}{{\ensuremath{\mathcal H}} }

\newcommand{\cD}{{\ensuremath{\mathcal D}} }

\newcommand{\bP}{{\ensuremath{\mathbf P}} }
\newcommand{\bE}{{\ensuremath{\mathbf E}} }

\DeclareMathSymbol{\leqslant}{\mathalpha}{AMSa}{"36} 
\DeclareMathSymbol{\geqslant}{\mathalpha}{AMSa}{"3E} 
\DeclareMathSymbol{\eset}{\mathalpha}{AMSb}{"3F}     
\newcommand{\dd}{\,\text{\rm d}}             
\newcommand{\sumtwo}[2]{\sum_{\substack{#1 \\ #2}}} 

\newcommand{\bbE}{{\ensuremath{\mathbb E}} }

\newcommand{\bbL}{{\ensuremath{\mathbb L}} }

\newcommand{\bbP}{{\ensuremath{\mathbb P}} }

\newcommand{\ga}{\alpha}

\newcommand{\gd}{\delta}
\newcommand{\gep}{\varepsilon}       

\newcommand{\go}{\omega}
\newcommand{\gto}{{\tilde\omega}}
\newcommand{\gO}{\Omega}
\newcommand{\gl}{\lambda}

\newcommand{\gs}{\sigma}

\makeatletter
\def\captionfont@{\footnotesize}
\def\captionheadfont@{\scshape}

\long\def\@makecaption#1#2{%
  \vspace{2mm}
  \setbox\@tempboxa\vbox{\color@setgroup
    \advance\hsize-6pc\noindent
    \captionfont@\captionheadfont@#1\@xp\@ifnotempty\@xp
        {\@cdr#2\@nil}{.\captionfont@\upshape\enspace#2}%
    \unskip\kern-6pc\par
    \global\setbox\@ne\lastbox\color@endgroup}%
  \ifhbox\@ne 
    \setbox\@ne\hbox{\unhbox\@ne\unskip\unskip\unpenalty\unkern}%
  \fi
  \ifdim\wd\@tempboxa=\z@ 
    \setbox\@ne\hbox to\columnwidth{\hss\kern-6pc\box\@ne\hss}%
  \else 
    \setbox\@ne\vbox{\unvbox\@tempboxa\parskip\z@skip
        \noindent\unhbox\@ne\advance\hsize-6pc\par}%
\fi
  \ifnum\@tempcnta<64 
    \addvspace\abovecaptionskip
    \moveright 3pc\box\@ne
  \else 
    \moveright 3pc\box\@ne
    \nobreak
    \vskip\belowcaptionskip
  \fi
\relax
}
\makeatother
\def\writefig#1 #2 #3 {\rlap{\kern #1 truecm
\raise #2 truecm \hbox{#3}}}

\newcommand{\spC}{\textsc{c}}
\newcommand{\tf}{\textsc{f}}

\newcommand{\bv}{{\underline{v}}}
\newcommand{\bw}{{\underline{w}}}
\newcommand{\bgo}{{\underline{\go}}}

\newcommand{\Ham}{{\cH}}
\newcommand{\Hamc}{{\Ham_{N,\bgo}(S)}}

\newcommand{\SRW}{S^{\textsf RW}}

\begin{document}
\date{October 4, 2005; accepted January 17, 2006}

\title[Disordered copolymers with adsorption]{The localized phase
of disordered copolymers with adsorption}

\author{Giambattista Giacomin}
\address{Laboratoire de Probabilit{\'e}s de P 6\ \& 7 (CNRS U.M.R. 7599)
  and  Universit{\'e} Paris 7 -- Denis Diderot,
U.F.R.                Math\'ematiques, Case 7012,
                2 place Jussieu 75251 Paris cedex 05, France}
\email{giacomin\@@math.jussieu.fr}
\urladdr{http://www.proba.jussieu.fr/pageperso/giacomin/GBpage.html}

\author{Fabio Lucio Toninelli}
\address{
Laboratoire de Physique, UMR-CNRS 5672, ENS Lyon, 46 All\'ee d'Italie, 
69364 Lyon Cedex 07, France}
\email{fltonine@ens-lyon.fr}
\urladdr{http://perso.ens-lyon.fr/fabio-lucio.toninelli}

\begin{abstract}
We analyze the localized phase of a general model
of a directed polymer in the proximity of an interface 
that separates two solvents. 
Each monomer unit carries a charge, $\go_n$, that determines the type 
(attractive or repulsive) and the strength
of its interaction with the solvents. 
In addition, there is a polymer--interface interaction  and
we want to model  the case in which there are impurities $\gto_n$,
called again charges, at the interface.
The charges are distributed in an inhomogeneous fashion
along the chain and at the interface: more precisely
the model we consider is of quenched disorder type.
\\
It is well known that such a model undergoes a localization/delocalization
transition. We focus on the localized phase, where the polymer sticks  to the interface.
Our new results include estimates on the exponential decay of correlations,
the proof that the free energy is infinitely differentiable
away from the transition and estimates on finite--size corrections to the 
thermodynamic limit of the free energy per unit site.
Other results we prove, instead, generalize earlier works that typically
deal either with the case of copolymers near an homogeneous interface ($\gto\equiv 0$) 
or with the case of disordered pinning, where 
the only  polymer--environment interaction is at the interface ($\go\equiv0$).
Moreover, with respect to most of the previous literature, we 
work with rather general distributions of charges
(we will assume only a suitable concentration inequality) and we allow 
more freedom on the law of the underlying random walk.
\end{abstract}

\maketitle

\section{Introduction}
\subsection{Copolymers, selective solvents and  adsorption}
\label{sec:premodel}
Polymers are repetitive chains of elementary blocks called monomers (or monomer units).
Copolymers 
are inhomogeneous polymers, in the sense
that each monomer unit carries
 a charge, and the charge is distributed along the chain 
in a {\sl disordered} way. It is well known that when
 the medium surrounding
 the copolymer  is made  of two solvents, separated
for example by an interface, and the solvents interact 
with the monomers according to the value of the charge,
the typical behavior of the copolymer may differ substantially
from the case in which the medium is homogeneous. On copolymers
there is an extremely extended literature, given above all
their practical relevance, see for example   \cite{cf:GHLO} and \cite{cf:G}
and references therein.
Moreover, for a realistic model of the interface, one should consider
the possibility of the presence of impurities 
or fluctuations in the interface layer, as in \cite{cf:SW}. 
As an extreme, but very important
example, one could consider also the case in which the interactions
at the interface are essentially the only relevant ones, see \cite{cf:DHV,cf:FLNO}.

In order to be more concrete, let us introduce a specific model,
which is just a particular example of the general class we 
consider. It is based on the process
 $\SRW=\{\SRW_n \}_{n=0,1,\ldots}$, a simple random walk on $\Z$,
with  $\SRW_0=0$ and  $\{ \SRW_j- \SRW_{j-1}\}_{j\in\N}$, 
$\N:=\{1,2, \ldots\}$,  a sequence
of IID random variables with  $\bP \left( \SRW_1= 1\right)=\bP \left( \SRW_1= -1\right)=
1/2$.

The process
$\SRW$ has to be interpreted as a directed polymer in $(1+1)$--dimension, and $\bP$ as its law in absence of
any interaction with the environment ({\sl free polymer}).
The polymer--environment interaction depends on the
{\sl charges} $\bgo =(\go, \tilde \go)\in \R ^{\N} \times
\R ^{\N} $ and on four real parameters
$\bv := (\gl, h, \tilde \gl, \tilde h)$: without loss of generality
we will assume $\gl$, $h$ and $\tilde \gl$ to be non--negative.
Let us set $S_n= \SRW_{2n}/2$ and  let us introduce the family of Boltzmann measures indexed by $N \in \N$ 
\begin{multline}
\label{eq:Boltzmann}
\frac {\dd \bP^{\bv} _{N, \bgo}} {\dd \bP} (\SRW)
\, =
\\
\frac 1
{\tilde Z_{N,\bgo}^\bv}
{\exp\left( \gl \sum_{n=1}^N \left( \go_n +h\right) \sign \left( S_n \right)+
\tilde \gl\sum_{n=1}^N \left( \gto_n +\tilde h\right)\ind_{\{S_n=0\}}
\right)} \ind_{\left\{ S_N=0\right\}},
\end{multline}
with the convention that $\sign \left(S_{n}\right)=\sign \left(\SRW_{2n-1}\right)$
for any $n$ such that $S_{n}=0$.
The superscript $\bv$ will be often omitted.

The model is completely defined once 
 $\bgo$ is given: being interested in the disordered case, we
 choose for instance    $\bgo$ an realization of an IID family,
 of law $\bbP$, of symmetric variables taking value $\pm 1$. 

\begin{figure}[!h]
\begin{center}
\leavevmode
\epsfysize = 7cm
\psfragscanon
\psfrag{n}[c][l]{$n$}
\psfrag{N}[c][l]{$N$}
\psfrag{0}[c][l]{$0$}
\psfrag{+}[c][l]{$+$}
\psfrag{-}[c][l]{$-$}
\psfrag{to1}[c][l]{$\tilde\omega_1$}
\psfrag{to2}[c][l]{$\tilde\omega_2$}
\psfrag{o1}[c][l]{$\omega_1$}
\psfrag{o2}[c][l]{$\omega_2$}
\psfrag{Sn}[c][l]{$S$}
\psfrag{Water}[c][l]{{\bf Water}}
\psfrag{Oil}[c][l]{{\bf Oil}}
\psfrag{n=1}[c][l]{$1$}
\psfrag{n=2}[c][l]{$2$}
\psfrag{n=3}[c][l]{$3$}
\psfrag{n=4}[c][l]{$4$}
\psfrag{n=5}[c][l]{$5$}
\epsfbox{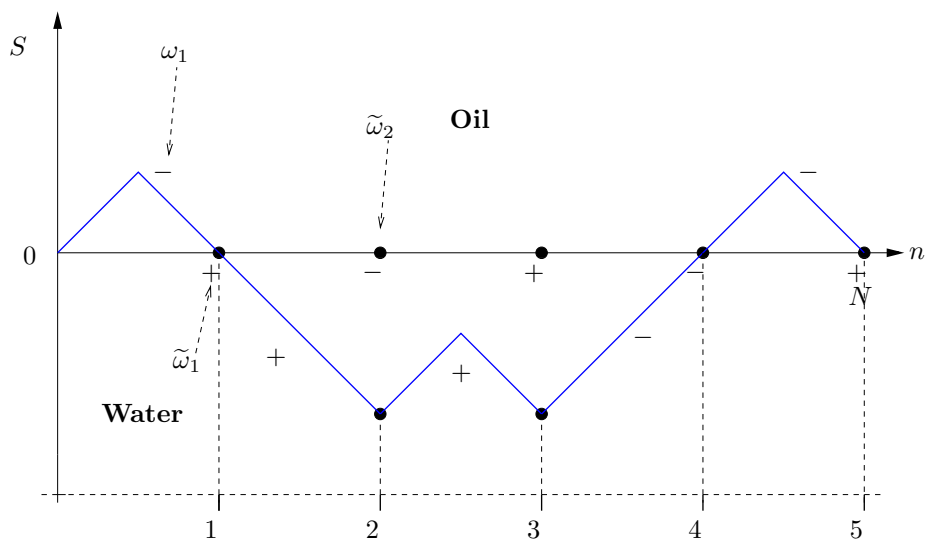}
\end{center}
\caption{\label{fig:modello} 
The polymer is at the interface between two solvents, say oil and water, situated
in the positive and negative half-planes, respectively. Along the (one--dimensional)
interface $S=0$ are placed random charges $\tilde\omega_n$, 
positive or negative. 
The polymer prefers  to touch the interface whenever
$\gto_n+\tilde h>0$. On the other hand, a random charge $\omega_n$ 
is associated to the $n^{th}$ monomer of the chain, i.e., to the 
portion of the polymer contained between $n-1$ and $n$. This is responsible for the polymer-solvent interaction:
if $\omega_n+h>0$ the $n^{th}$ monomer prefers to be in oil, otherwise it prefers to be in water. 
The non--negative parameters $\lambda,\tilde\gl$ may be thought of as two effective
inverse temperatures, and $h$ as a measure of the asymmetry between the two solvents.
Note that by convention, if $S_n=0$ then $\sign(S_n)$ coincides with the sign of $S_{n-(1/2)}$, which
is unambiguously defined. For instance,
$\sign(S_1)=+1$ and $\sign(S_4)=-1$.
In this picture, $N=5$ and in fact the endpoint $S_N$ is pinned at zero. 
}
\end{figure}

We invite the reader to look at Figure \ref{fig:modello} in order to get an intuitive idea on this model.

\smallskip

\subsection{The model}
\label{sec:model}
An important observation on the model we have introduced
is that its Hamiltonian may be easily rewritten in terms
of the sequence $\tau :=\{\tau_i\}_i$, defined by setting 
 $\tau_0:=0$ and, for $i\in \N$,    
$\tau_{i}:= \inf\left\{n> \tau_{i-1}:\, S_n=0 \right\}$ 
($\tau_{i} < +\infty$ with probability one for every $i$ since $\SRW$ is recurrent)
 and in terms of the sign of the excursions, that, 
 conditionally on $\tau$, are just an   independent  sequence of IID symmetric
 variables taking values $\pm 1$.
 Of course 
 the (strong) Markov property of  $\SRW$ immediately yields also that 
$\{ \tau_j-\tau_{j-1}\}_{j\in \N}$ is an IID sequence.
\medskip

We are now going to introduce  the general class of models that we 
consider. These models are based on 
a real--valued free process $S:= \left\{S_n\right\}_{n=0,1,\ldots}$, with law $\bP$  
that satisfies the following properties: 
\smallskip

\begin{enumerate}
\item The sequence $\tau=\{\tau_j\}_{j=0,1, \ldots}$ of successive returns to $0$ 
is an infinite sequence with $\tau_0=0$, so $S$ is a process starting from $0$
and for which $0$ is a recurrent state. Moreover $\tau$ is a renewal sequence,
that is $\{ \tau_{j}-\tau_{j-1}\}_{j=1,2, \ldots}$ is a sequence of IID random variables,
and we set $K(n):=\bP(\tau_1=n)$. 
\item 
For some 
 $\ga \ge 1$ 
and for some function 
$L(\cdot)$  which is  slowly varying at infinity (see below for
the definition and properties of slowly varying functions)
 \begin{equation}
 \label{eq:alpha}
  K(n)\, =\, \frac{L(n)}{n^\alpha}.
 \end{equation} 
In particular  $K(n)>0$ for every $n\in \N$. 
\item
For $i\in \N$ such that $\tau_i -\tau_{i-1}>1$ we set $s_i:= \sign(S_{\tau_i -1})$.
Then conditionally on $\tau$, $\{s_i\}_{i: \, \tau_i -\tau_{i-1}>1}$ is an IID sequence
of symmetric random variables taking values $\pm 1$.
Conventionally we complete the sequence $\{s_i\}_{i\in \N}$ (i.e., we choose the $s_i$ for $i$ such that 
$\tau_{i+1}-\tau_i=1$) by tossing  independent
fair coins and, always by convention, we stipulate that $\sign(S_{\tau_i })=s_i$.
Of course, conditionally on $\tau$,  $\{s_i\}_{i\in \N}$ is just independent fair coin tossing.
\item Conditionally on $\tau$, $\{ (S_{\tau_{i}+1}, \ldots, S_{\tau_{i+1}}) \}_{i=0,1,\ldots}$
is an independent sequence of random vectors. Moreover the law of 
$(S_{\tau_i+1}, \ldots, S_{\tau_{i+1}})$, conditionally on $\tau$,
depends  on $\tau$ only via the value of $\tau_{i+1}-\tau_i$.
Note that this property implies that, conditionally on $S_m=0$, 
the process $\{S_n\}_{n\ge m}$ has the same distribution as 
the original process $\{S_n\}_{n\ge0}$. With some abuse of language, we will
call this property the {\em renewal property of $S$}.
 \end{enumerate}
\medskip

The free process $S$ is put  in interaction with an environment
via  the 
{\sl charges} $\bgo =(\go, \tilde \go)\in \R ^{\N} \times
\R ^{\N} $.
The definition of 
$ \bP^{\bv} _{N, \bgo}$ is still
as in \eqref{eq:Boltzmann}, provided  one suppresses 
the superscript RW in the left--hand side.
The process $S$ constructed in Section 
\ref{sec:premodel} starting from $\SRW$ corresponds to the case
of $\ga=3/2$ and $\lim_{n\to \infty}L(n)=1/(2\sqrt{\pi})$.

We recall that a function $L(\cdot) $ is slowly varying (at infinity) 
if   $L(\cdot)$ is a measurable function from
$(0, \infty)$ to $ (0, \infty)$ such that
 $\lim_{r\to \infty} 
L(x r)/L(r)=1$ for every $x>0$. 
One of the properties of
slowly varying functions is that  both $L(r)$ and $1/L(r)$
are $o(r^\epsilon)$ for every $\epsilon >0$. An example of slowly varying function
is 
 $r \mapsto (\log (r+1))^b$, for  $b\in \R$, 
 but also $r \mapsto \exp( (\log (r+1))^b)$,  for $b <1$, 
as well as any positive function for which $\lim_{r\to \infty} L(r)>0$. 
A complete treatment of  slowly varying functions is 
found in
 \cite{cf:Feller2}, but these  functions are needed
for us to in order to work in  a reasonably general and well defined set--up:
 we will use no fine property of slowly varying function
since in most of the cases rough bounds on $K(\cdot)$
will suffice. 

\medskip
\begin{remark}
\label{rem:S}
\rm
In short, one may think of building $S$
by first assigning the return times to zero according to a renewal process.
The excursions are then {\sl glued} to the renewal points, 
but essentially the only relevant 
aspect of the excursions for us is the sign, that is chosen by repeated tossing of a fair coin.
Note again that the energy of the model depends only on $\tau$ and
on $\{s_i\}_i$, and  not on the details of the excursions in the upper or lower half plane.
The last property in the list above, the renewal property of $S$, 
 makes a bit more precise  
what the excursions of the process really are: this property 
is very useful to have a nice pictorial vision of the process, but it
is rather inessential for us. We will use it    only 
in stating Theorem~\ref{th:correlazioni} since this way it turns out to be 
somewhat nicer and more intuitive, but the essence of
our analysis lies in $\tau$. 
Note also that if $\gl=0$ the energy does not depend on $\{s_i\}_i$.
In this case we are dealing with a pure pinning model and 
insisting on $S$ taking values in $\R$ is a useless limitation. 
\end{remark}
 \medskip

\begin{remark}
\label{rem:S1}
\rm
We have chosen $0$ to be a recurrent state of $S$ because it simplifies 
a little the notations, but everything carries over to the  case in which 
$0$ is transient. 
\end{remark}
 \medskip

The sequence $\go$ is chosen as a typical realization of an IID
sequence of random variables, still denoted by $\go = \{ \go_n \}_n$.  
We make the very same assumptions on 
 $\gto$ and, in addition,  $\omega$ and $\gto$ are independent. The law of $\bgo$ will be denoted by $\bbP$, and the 
corresponding expectation by $\bbE$.
A further assumption on $\bgo$ is the following {\sl concentration inequality}:
there exists
 a positive constant $C$ such that for every $N$,
 for every Lipschitz and convex function $g: \R ^{2N} \to \R$ with 
 $g(\bgo) := g\left(\go_1, \ldots, \go _N,\gto_1,\ldots,
\gto_N \right) 
\in  L^1 \left( \bbP\right)$ and  for $t\ge 0$
 \begin{equation}
 \label{eq:Lip}
 \bbP
 \left(  \big \vert g(\bgo) - \bbE\left[ g(\bgo)\right]\big\vert \ge t\right)
 \, \le \,  C \exp \left( -\frac{t^2}{C \Vert g \Vert ^2 _{\rm Lip}}\right),
 \end{equation}
 where
  $\Vert g \Vert_{\rm Lip}$ is the Lipschitz constant of $g$ with respect to the Euclidean distance on $\R^{2N}$.

Of course such an inequality implies that $\go_1$ and $\tilde \go_1$
are exponentially integrable: without loss of generality we   
assume $\go_1$ and $\tilde \go _1$ to be centered and of unit variance.

  \medskip
  
\begin{remark}\rm
In practice, it is sufficient to check that the concentration inequality holds for $\go$ and for $\gto$ separately.
In fact, suppose that 
\begin{equation}
 \label{eq:Lip2}
 \bbP
 \left(  \big\vert G(\go) - \bbE\left[ G(\go)\right]\big\vert \ge t\right)
 \, \le \,  C \exp \left( -\frac{t^2}{C \Vert G \Vert ^2 _{\rm Lip}}\right),
 \end{equation}
for every $G: \R ^{N} \to \R$ and $G(\go):=G(\go_1,\ldots,\go_N)$, and similarly for $\gto$.
Then, 
\begin{multline}
  \label{eq:trick}
  \bbP\left(\big\vert g(\bgo)-\bbE[g(\bgo)]\big\vert \ge u\right)\le 
  \\
\bbP\left(\big\vert \bbE[g(\bgo)\vert \go]-\bbE[g(\bgo)]\big\vert \ge \frac u2\right)
+  \bbP\left(\big\vert g(\bgo)-\bbE[g(\bgo)\vert \go]\big\vert \ge \frac u2\right).
\end{multline}
Using twice  \eqref{eq:Lip2}, once for $G_1(\go):=\bbE [g(\bgo)\vert \go]$ and once for
$G^{\go}_2(\gto):=g(\go,\gto)$, noting that $\max(||G_1||_{\rm Lip} 
, ||G^{\go}_2||_{\rm Lip})\le ||g||_{\rm Lip}$ 
uniformly in $\go$, one obtains  immediately \eqref{eq:Lip}.

The concentration inequality \eqref{eq:Lip} is known to hold with
a certain generality: its validity for
 the Gaussian case and  
for the case of bounded random variables  is by now 
a classical result (\cite{cf:LedouxAMS}). While of course such an inequality requires
  $\bbE \left[ \exp\left( \gep \left( \go_1 ^2 + \gto_1^2\right) \right) \right] <\infty$
for some $\gep>0$, a complete characterization of the distributions for which  
 \eqref{eq:Lip} holds is, to our knowledge, still lacking, 
 but among these distributions there are,
 for example,
the cases in which the laws of $\go_1$ and  $\gto_1$ satisfy the log--Sobolev inequality, see \cite{cf:newlook}, 
\cite{cf:Ledoux} and \cite{cf:cedric}, therefore, in particular, whenever 
 $\go_1$ and  $\gto_1$ are continuous variables with positive densities of the
 form $\exp(-V(\cdot))$, with  
 $V \in C^2$ and $V''(\cdot)\ge c> 0$ when restricted to $ (-\infty, -K)\cup 
 (K , \infty)$, for some $K>0$.
\end{remark}

\medskip

\begin{remark}
\rm
We have chosen to work assuming concentration because it provides 
a unified rather general framework in which proofs are at several instances  
much shorter (and, possibly, more transparent).
However, as it will be clear,
several results hold under weaker assumptions.
On the other hand we  deal only with 
the polymer pinned at the endpoint, that is, constrained to $S_N=0$:
we have chosen this case for the sake of conciseness, but we could have
decided for example to  leave the endpoint free. 
\end{remark}

\medskip

\subsection{Free energy, localization and delocalization}
\label{subs:fe}
Under the above assumptions on the disorder
the {\sl quenched free energy} of the system exists, namely
the limit
\begin{equation}
\label{eq:free_energy}
f(\bv)\, := \, \lim_{N\to \infty} \frac 1N \log  \tilde Z^\bv_{N,\bgo},
\end{equation}
exists  $\bbP (\dd \bgo)$--almost
surely and in the $\bbL^1\left( \bbP\right)$ sense.  
The existence of this limit
can be proven via standard super--additivity arguments (we refer for example to \cite{cf:G} for details).
We stress that the concentration inequality implies immediately
that $f(\bv)$ is not random, but such a result may be proven under much
weaker assumptions, see e.g. \cite{cf:G}. 
\smallskip
 
 A simple but fundamental observation is that
\begin{equation}
\label{eq:delocfe}
f(\bv)\, \ge \, \gl h.
\end{equation} 
The proof of this
 is elementary:
 if we set $\gO_N^+= \{S\in \Omega:\, S_n>0$ for $ n=1, 2, \ldots , N-1$ and  $
 S_N=0\}$ 
 we have
\begin{multline}
\label{eq:step_deloc}
\frac 1N \bbE \log  \tilde Z^\bv_{N,\bgo}  \ge 
\\ 
\frac 1N \bbE \log \bE 
\left[ 
\exp\left( \gl \sum_{n=1}^N \left( \go_n +h\right) \sign \left( S_n \right)
+\tilde \gl\sum_{n=1}^N \left(\gto_n +\tilde h\right)\ind_{\{ S_n =0 \} }
\right)
;  \gO_N^+
\right]
\\
=\,
\frac \gl N {  \sum_{n=1}^N \left(\bbE\left[\go_1 \right]+h\right) } 
 +
\frac {\tilde \gl}N \left( \bbE \left[\tilde \go _N \right] +\tilde h \right)
+
\frac 1N 
\log \bP \left( \gO_N^+\right)\, \stackrel{N \to \infty}{\longrightarrow}\, \gl h,
\end{multline}
where 
we have  applied the
 fact that, by  \eqref{eq:alpha}, 
 $\log \bP\left( \gO_N^+\right)=\log((1/2) K(N))=o(N)$, for $N\to\infty$.

The observation \eqref{eq:delocfe}, above all if viewed in the light
of its proof, suggests the definition $\tf (\bv):=f(\bv)- \gl h$ and 
the following
partition of the parameter space (or {\sl phase diagram}):

\smallskip
\begin{itemize}
\item The localized region: $\cL = \left\{ \bv : \, \tf(\bv)> 0\right\}$;
\item The delocalized region:  $\cD = \left\{ \bv: \, \tf(\bv)=0\right\}$.
\end{itemize}
\smallskip

Along with this definition we observe that 
\begin{equation}
\label{eq:Boltzmann1}
\frac {\dd \bP^{\bv} _{N, \bgo}} {\dd \bP} (S)
\, \propto \, 
\exp\left(-2 \gl \sum_{n=1}^N \left( \go_n +h\right) \Delta_n+
\tilde \gl\sum_{n=1}^N \left( \gto_n +\tilde h\right)\delta_n
\right)\delta_N,
\end{equation}
where 
\begin{equation}
\Delta_n := \ind_{\{\sign(S_n)=-1\}}= \left(1-\sign(S_n)\right)/2 \ \ \ \text{ and }
\ \ \
 \delta_n:=
\ind_{\{S_n=0\}}.
\end{equation}
Of course the normalization constant 
\begin{equation}
\label{eq:normalisaz}
Z_{N,\bgo}^\bv \, :=\, \bE\left[ \exp\left(-2 \gl \sum_{n=1}^N \left( \go_n +h\right) \Delta_n+
\tilde \gl\sum_{n=1}^N \left( \gto_n +\tilde h\right)\delta_n
\right)\delta_N\right],
\end{equation}
changes, but $Z_{N,\bgo}^\bv =\exp\left(  -\gl \sum_{n=1}^N (\go_n + h)\right)
\tilde Z_{N,\bgo}^\bv$ so that the $\bbP(\dd \bgo)$--a.s. and $\bbL^1(\bbP)$ asymptotic behavior of
$(1/N)\log Z_{N,\bgo}^\bv $ are given by $\tf (\bv)$.
We will always work with $Z_{N,\bgo}^\bv $, 
in order to conform with most of the previous mathematical literature.

\subsection{Discussion of the model}
\label{sec:discuss}
For the model we introduced there is a very vast  literature, mostly in chemistry,
physics and bio--physics. 
It often goes under the name of  {\sl copolymer  with adsorption}, see e.g.
\cite{cf:SW} and references therein, and  
such a  name clearly reflects the superposition of  two distinct 
polymer--environment 
interactions:
\smallskip

\begin{itemize}
\item The monomer--solvent interaction, associated to the charges $\go$. 
Some mon\-o\-mers prefer
one solvent and some prefer the other one. Since the charges are placed in an 
inhomogeneous 
way along the chain, energetically favored trajectories
need to stick close to the interface. Whether pinning actually takes place or not depends on the interplay between
energetic gain and entropic loss associated to localization 
(trajectories that stay close to the interface have a much smaller entropy
than those which wander away).

If $\gl>0$ and $\tilde \gl =0$ only this interaction is present and we will
call the model simply {\sl copolymer}. 

\item The monomer--interface interaction, associated  to  $\gto$. This interaction
leads  to a pinning (or depinning)
phenomenon with a more direct mechanism:  trajectories 
are energetically favored
if they touch  the interface {\sl as often as possible}
at points where $\gto_n+\tilde h> 0$,
avoiding at the same time the points 
in which $\gto_n+\tilde h< 0$. Also in this case, 
a non--trivial energy/entropy competition is
responsible for the localization/delocalization transition.

 When $\gl=0$ and $\tilde \gl>0$ we will refer to the model as {\sl pinning} model. 
\end{itemize} 

\smallskip


The copolymer model  has received a lot of attention:
we mention in particular \cite{cf:GHLO}, in which it
 was first introduced and the replica method was applied in order to investigate the transition.
Rigorous work started with \cite{cf:Sinai},
followed by \cite{cf:AZ}: these works deal with the case
$h=0$ and $\go_1$ symmetric and taking only the values $\pm 1$
({\sl binary charges}).
It turns out that in such a case there is no transition and the model
is localized for every $\gl >0$.

In general,  one may distinguish between results concerning
the free energy and results on pathwise behavior of the polymer.
About the first point, we mention that in \cite{cf:BdH}  the model with $h \ge 0$  has been considered,
still with the  choice of binary charges, and
the existence of a transition has been established, along with
estimates on the critical curve and remarkable
limiting properties of the free energy of the model in
the limit of weak coupling ($\gl $ small). 
Improved estimates on free energy and critical curve may be found in \cite{cf:BG}.
In the physical literature one can find 
a number of conjectures, mostly on the free energy behavior,
that are far from clarifying  the phase diagram. 
In this respect, it is interesting to mention that  recent 
numerical simulations (\cite{cf:CGG}) show that the critical line is different from that 
predicted in the theoretical physics literature, which means that the 
localization mechanism is still poorly understood.

\smallskip

Disordered pinning has been extensively studied in the 
physical literature, see e.g. \cite{cf:DHV} and \cite{cf:FLNO}
(see also \cite{cf:GT05} for more recent references),
but much less in the mathematical one. However the
model has started attracting attention lately,
see  \cite{cf:AS}, \cite{cf:Petrelis} and \cite{cf:GT05}.
We should stress that there is no agreement 
in the physical literature on several important issues
for disordered pinning. For instance, it is still unclear whether the critical curve coincides with the
so--called annealed curve.

\smallskip

About the study of path behavior, 
there is a basic difference between the localized and the delocalized phase: in the 
first case, since $\tf(\bv)>0$, the interaction produces an exponential modification of the free polymer
measure, and therefore 
{\sl Large Deviation} techniques apply very naturally. The path behavior of the  
copolymer model in the localized phase
  has been 
considered in
\cite{cf:Sinai,cf:AZ}, for $h=0$ and binary charges, while  in 
\cite{cf:BisdH} also the case $h>0$ is taken into account.
In \cite{cf:BisdH} the focus is on the
  Gibbsian characterization of the infinite volume polymer  measure
(in the localized phase). The delocalized phase is more subtle, due to the fact
that $\tf(\bv)=0$, and path delocalization estimates involve estimates on {\sl Moderate Deviations} of the
free energy. Results on
the path behavior in the delocalized phase have been obtained only recently
in \cite{cf:GT}, both for the copolymer and for the disordered pinning model.

\section{Main results}

In the present work, we consider the localized phase of the general model
defined in Section~\ref{sec:model} and formula
 \eqref{eq:Boltzmann}. 
In addition to giving new results, our approach provides  also a setting to reinterpret in a simpler way 
known results for copolymer and pinning models.

\subsection{Smoothness of the free energy and decay of correlations}
The free energy is everywhere continuous and almost everywhere differentiable, by convexity,
so in particular $\cL$ is an open set. However,
one can go much beyond that, as our first theorem shows:
\medskip

\begin{theorem}
\label{th:cinf}
   $\tf$ is infinitely differentiable in $\cL$.
\end{theorem}
\medskip

An interesting problem is to study the regularity properties of $\tf(\cdot)$ at the boundary between $\cL$ and 
$\cD$, where it is  non--analytic. This corresponds to investigating
 the order of the
localization/delocalization transition. Recently, an important step in this direction was performed in
\cite{cf:GT05} and \cite{cf:GTlett} where it was proved, in particular,  that the first derivatives of $\tf(\cdot)$ 
are continuous on the boundary.
In other words, the (de)localization transition is at least of second order.

\smallskip

As it will become clear in Section \ref{sec:reg}, the smoothness of the free energy in $\mathcal L$ boils down to  
a property of exponential decay of (average) correlations. 
 For this implication we essentially rely on  
 \cite{cf:griffiths}, where a similar result has been 
 proven in the context of disordered Ising models.

Let us therefore state the  decay of correlation property.
We say that  $A$ is a bounded local observable if   
 $A$ is a real bounded measurable local function of the path configurations.
 In the sequel $\mathcal S(A)$ will denote the {\sl support of} $A$, that is
 the intersection of all the  subsets $I\subset \N$ of the form $I=\{\ell,\ldots,k\}$, with $k,\ell\in \N$, 
 such that $A$ is measurable with respect to the $\gs$--algebra  $\sigma(S_n:\, n\in I)$.
 \smallskip

  We have
  \bigskip
  
\begin{theorem}
\label{th:correlazioni}

For every 
$\bv\in \cL$ there exist finite constants $c_1(\bv),c_2(\bv)> 0$ such that 
the following holds for every $N\in \N$: 
\begin{itemize}
\item  (Exponential decay of correlations.) For every couple of bounded local observables $A$ and $B$ we have 
\begin{multline}
\label{eq:decad1}   
\phantom{movemo} 
\bbE \left[\left|\bE_{N,\bgo}(A B)-\bE_{N,\bgo}(A)\bE_{N,\bgo}(B)
\right|\right]
\\
\le\,  c_1 \Vert A\Vert_\infty\Vert B\Vert_\infty
\exp(-c_2 d(\mathcal S(A),\mathcal S(B)))
  \end{multline}
where $d(I,J)=\min\{|i-j|, i\in I,j\in J\}$ if $I,J\subset \N$.
\item (Influence of the boundary.)
For every bounded local observable $A$ and $k\in \N$, such that
$\mathcal S(A)\subset \{1,\ldots,k\}$, we have 
\begin{equation}
\label{eq:bordo}
\sup_{N>k}
  \bbE\left[\left| \bE_{N,\bgo}(A)-\bE_{k,\bgo}(A)
\right|\right]\le  c_1 \Vert A\Vert_\infty
\exp(-c_2d(\mathcal S(A),\{k\})).
\end{equation}
\item  For every bounded local observable $A$ 
the following limit exists $\bbP(\dd \bgo)$--almost surely:
  \begin{equation}
    \label{eq:exlim}
\lim_{N\to\infty} \;  \bE_{N,\bgo}(A)\, =: \,  \bE_{\infty,\bgo}(A).
  \end{equation}
\end{itemize}
\end{theorem}
\bigskip

\begin{remark}
\label{rem:a.s.}
\rm
From \eqref{eq:decad1}
one may easily extract 
an almost sure statement. Choose two bounded local observables $A$
and $B$ and 
 set  $ B_k (S)= B(\theta ^k S)$, 
i.e. $(\theta S)_n= S_{n+1}$.
Take the limit  $N\to \infty $ in \eqref{eq:decad1} to obtain
\begin{equation}
\label{eq:decad1foras}    
\bbE \left[\exp(c_2 k/2)\left|\bE_{\infty,\bgo}(A B_k)-\bE_{\infty,\bgo}(A)\bE_{\infty,\bgo}(B_k)
\right|\right]\le 
\exp(-c_2 k/4),
  \end{equation}
for $k$ sufficiently large.   Therefore the Fubini--Tonelli theorem 
and \eqref{eq:decad1foras}   yield
\begin{equation}
\bbE \left[ \sum_k
\exp(c_2 k/2)\left|\bE_{\infty,\bgo}(A B_k)-\bE_{\infty,\bgo}(A)\bE_{\infty,\bgo}(B_k)
\right|
\right]\, <\, \infty.
\end{equation}
The series appearing in the left--hand side is therefore $\bbP(\dd \bgo)$--a.s. convergent. This implies that
there exists a random variable $C_{A,B}(\bgo)$, 
$C_{A,B}( \bgo)< \infty$ $\bbP(\dd \bgo)$--a.s., such that
\begin{equation}
\label{eq:a.s.stat}
\left|\bE_{\infty,\bgo}(A B_k)-\bE_{\infty,\bgo}(A)\bE_{\infty,\bgo}(B_k)
\right| \le C_{A,B}(\bgo) \exp(-c_2 k/2),
\end{equation}
for every $k$.
\end{remark}

\subsection{Path localization and maximal excursions}
\label{sec:path}

We consider now the question of whether or not
knowing that $\bv \in \cL$ does mean that the path of the polymer 
is really tight to the interface. Even if this question has not
been treated for the general model  we are considering here,
the techniques used in \cite{cf:Sinai}, \cite{cf:AZ} and \cite{cf:BisdH}, see also 
\cite{cf:G}, 
may be applied directly and one obtains for example that, in the case $S_n=\SRW_{2n}/2$,
for every $\bv \in \cL$
  there exist finite constants $c_1,c_2> 0$ such that for every $N$, $s\in \N$, and $0\le k\le N$
  \begin{equation}
\label{eq:altezza}
    \bbE\,\bP_{N,\bgo}\left(|S_k|\ge s\right)\le c_1 e^{-c_2 s},
  \end{equation}
  or one can obtain an analogous $\bbP(\dd \bgo)$--a.s. result, which is 
  a bit more involved to state  \cite{cf:Sinai}.

  \smallskip
  
The reason for  revisiting this type of results, besides generalizing them to our case,
is that they are only bounds and we would like to find estimates that
are sharp to leading order.
A notable exception 
is the case of some of the results in \cite[Th. 5.3 and Th. 6.1]{cf:AZ} where
the precise asymptotic size of the largest excursion (and of the maximum 
displacement of the chain from the interface) is obtained. 
This result is a bit surprising since it depends on a certain  annealed decay exponent.
This exponent turns out to be
 different from the decay exponent one finds for $\bbP(\dd \bgo)$--a.s. estimates, see discussion after our
Theorem \ref{th:zu}.
The argument of the crucial point of the proof
of  \cite[Th. 5.3]{cf:AZ} looks obscure to us and we propose here a different one, based on decay
of correlations. 

\smallskip

Of course we present our results in terms of {\sl excursion lengths}.
 Recall the definition of the return times $\tau$ in Section \ref{sec:model}. 
 For every $k\in \{1, \ldots, N-1\}$ let us set
 $\boldsymbol \Delta_N(k)=\tau_{\iota (k)+1}-\tau_{\iota (k)}$, with $\iota (k)$ equal 
 to the value $i$ such that
  $k \in \left\{ \tau_{i}, \ldots, \tau_{i+1}-1\right\}$.
  So $\iota(k)$ is the left margin of the excursion to which $k$ belongs and 
  $\boldsymbol \Delta_N(k)$ is the length of such an excursion.

Two distinct questions can be posed concerning polymer excursions: one may be interested 
about the typical length of a given excursion or, more globally, about the typical length of the {\sl longest}
excursion, $\boldsymbol \Delta_N :=\max _{k}\boldsymbol \Delta_N(k)$.
Denote by $\theta $ the left shift on $\go$, like for $S$: 
$\{\theta\go\}_k=\go_{k+1}$.
About the first problem, we can prove: 
\bigskip

\begin{proposition}
\label{th:as} Take $ \bv \in \cL$ and let $\bgo$ be the two--sided sequence of IID random variables,
$\bgo:=(\go,\gto)\in \R^{\Z}\times\R^{\Z}$, with law  $\bbP$.
For every $\gep>0$ there exist  random variables $C^{(j)}_\gep(\bgo)$,
 $j=1,2$, such that $\bbP( 0< C^{(j)}_\gep (\bgo)<\infty)=1$ and
\begin{equation}
\label{eq:c1c2}
\begin{split}
 \bP_{N, \bgo}
\left( \boldsymbol \Delta_N (k) = s\right) \, &\ge\,
C^{(1)}_\gep(\theta^k\bgo) 
\exp\left( -\left( \tf (\bv ) +\gep \right)s \right),
\\
\bP_{N, \bgo}
\left( \boldsymbol \Delta_N (k) = s\right)
\, &\le \, C^{(2)}_\gep(\theta^k\bgo) 
\exp\left( -\left( \tf (\bv ) -\gep \right)s \right),
\end{split}
\end{equation}
for every $N$, every $k \in \{1, \ldots\, N-1\}$ and every $s\in \N$: for the first inequality,
the lower bound, we require also $s<N/2$.
\end{proposition} 
\medskip
 
Note that, in the definition of $\bP_{N,\bgo}$, the fact that $\bgo$ is a doubly infinite sequence is 
completely irrelevant,
since  the polymer measure depends only on $\bgo_1,\ldots,\bgo_N$. The introduction of two--sided
disorder sequences, which might seem a bit unnatural,
is needed here to have the almost--sure result uniformly in $N,k$ and  $s$.
\medskip

About the maximal excursion, we have: 

\medskip

\begin{theorem}
\label{th:zu}
The following holds: 
\begin{enumerate}
\item the limit
\begin{equation}
\label{eq:mu}
-\lim_{N \to \infty} \frac 1N \log \bbE
\left[ \frac {1+e^{-2\lambda\sum_{n=1}^N(\go_n+h)}} {Z_{N, \bgo}^\bv} \right] \, :=\, \mu (\bv) ,
\end{equation}
exists and satisfies the bounds $0 \le \mu (\bv) \le \tf( \bv)$. 
Moreover,  for every $c_1>0$ there exists 
$c_2>0$ such that $ c_2\left(\tf(\bv)\right)^2\le \mu(\bv)$, for $\tf(\bv)\le 1$
and $\max(\gl , \tilde \gl) \le c_1$.
\item
Fix  $\bv \in \cL$. For every $\gep\in (0,1)$ we have that
\begin{equation}
\label{eq:as1b}
\lim_{N \to \infty}  \bP_{N, \bgo}^\bv \left(
 \frac{\boldsymbol \Delta_N}{\log N} <\frac{1-\gep} {\mu (\bv)} \right)
 \, =\, 0,  \ \ \  \bbP(\dd \bgo) \text{--a.s.},
\end{equation}
and 
\begin{equation}
\label{eq:problb}
\lim_{N \to \infty}  \bP_{N, \bgo}^\bv \left(
 \frac{\boldsymbol \Delta_N}{\log N} >\frac {1+\gep}{\mu (\bv)}\right) 
 \, =\, 0, \ \ \ \text{in probability}.
\end{equation}
\end{enumerate}
\end{theorem}

\bigskip

In the localized region,  
under rather general conditions on the law $\bbP$, we can prove that $\mu(\bv)< \tf(\bv)$, see Appendix
\ref{sec:mu}. Notice therefore the gap between \eqref{eq:as1b}
and the result in Proposition~\ref{th:as}, so that the largest excursion 
appears to be achieved in atypical regions.

\medskip

\begin{remark}\rm
When the law of $\omega_1$ is symmetric, one can also prove (see Appendix \ref{sec:mumu})
that $\mu(\bv)$ is equivalently given by
\begin{equation}
\label{eq:mu2}
\mu (\bv)=-\lim_{N \to \infty} \frac 1N \log \bbE
\left[ \frac 1 {Z_{N, \bgo}^\bv} \right].
\end{equation}
This coincides with the expression given in \cite{cf:AZ}, where only the case
$\tilde \lambda=0$, $h=0$
and $\go_1$ taking values in $\{-1,1\}$ 
 was
considered.
\end{remark}

\subsection{Finite--size corrections and  central limit theorem for the free energy}

In this section, we investigate finite--size corrections to the infinite volume limit of the 
quenched average of the free energy, and the behavior of its disorder fluctuations.
About the first point, it is quite easy to prove (see also \cite[Proposition 2.6]{cf:AZ})
\begin{proposition}
\label{th:Nfin1}
  There exists $c_1<\infty$ such that, for every $\bv$ and $N\in \N$, one has
  \begin{eqnarray}
    \label{eq:Nfin1}
    \left|\tf(\bv)-\frac1N \bbE \log Z^{\bv}_{N,\bgo}\right|\le c_1 \frac{\log N}N +
\frac{\tilde \lambda}N (1+|\tilde h|).
  \end{eqnarray}
\end{proposition}
This bound is  somehow optimal in general, 
in the sense that it is possible to prove a lower bound of order $\log N/N$ if $\bv$
is in the {\em annealed region}, i.e., the sub--region of $\cD$  where $(1/N)\log \bbE Z^{\bv}_{N,\bgo}\to 0$.
However, in the localized region we can go much farther:
\begin{theorem}
\label{th:Nfin2}
Assume that $\bv\in \cL$. Then, there exists $c(\bv)<\infty$ such that, for every  $N\in \N$, one has
  \begin{eqnarray}
    \label{eq:Nfin2}
    \left|\tf(\bv)-\frac1N \bbE \log Z^{\bv}_{N,\bgo}\right|\le \frac{c(\bv)}N.
  \end{eqnarray}
\end{theorem}

About fluctuations, in \cite{cf:AZ} it was proven that, for $\tilde \lambda=0$, $h=0$
and $\go_1$ taking values in $\{-1,1\}$, the free energy satisfies in the large volume limit
a central limit theorem on the scale $1/\sqrt N$. Here, we generalize this result to the entire localized region.
Our proof employs basically the same idea as  in 
 \cite{cf:AZ}; however, the use of concentration of measure ideas,
plus a more direct way to show that the limit variance is not degenerate, allow for remarkable simplifications.
The precise result is the following:
\medskip

\begin{theorem}
\label{th:CLT}
  If $\bv \in \cL$ the following limit in law holds: 
  \begin{eqnarray}
    \frac1{\sqrt N}\left ( \log Z_{N,\bgo}^\bv -
    \bbE\log Z_{N,\bgo}^\bv\right)
\stackrel{N\to \infty}\longrightarrow\mathcal N(0,\sigma^2),
  \end{eqnarray}
with $\sigma^2=\sigma^2(\bv)>0$.
\end{theorem}
\medskip

\subsection{Some notational conventions}
\label{sec:conv}
For compactness we introduce a notation for the Hamiltonian
\begin{equation}
\Hamc \, =\, -2 \gl \sum_{n=1}^N \left( \go_n +h\right) \Delta_n+
\tilde \gl\sum_{n=1}^N \left( \gto_n +\tilde h\right)\delta_n,
\end{equation}
and for $\overline{\gO} \in \gs \left( S_n:\, n\in \N\right)$
we set
\begin{equation}
Z_{N, \bgo}\left(\overline{\gO} \right)\, =\,
\bE\Big[ \exp\left(\Hamc\right) ; \,\left\{
 S_N=0\right\} \cap \overline{\gO} \Big].
\end{equation}
Moreover, we set 
\begin{equation}
  \label{eq:zeta}
  \zeta(x)=e^{\tilde\gl(x+\tilde h)}.
\end{equation}


\section{Decay of correlations}
\label{sec:corr}

The proof of Theorem \ref{th:correlazioni} is based on the following lemma, which is somehow similar in spirit
to Lemmas 4 and 5 in \cite{cf:BisdH}.
\begin{lemma}
  \label{lemma:ritorni}
For every $\bv\in \cL$, there exist  constants $0<c_1(\bv),c_2(\bv)<\infty$ such that,
for every $ N\in \N$, $1\le a<b\le N$, $k\le (b-a-1)$ and $a< i_1<i_2\ldots <i_k< b$,
 \begin{eqnarray}
\label{eq:ritorni}
   \bbE\,\bP_{N,\bgo}(\{a \le j\le  b : S_j=0\}=\{a,i_1,i_2,\ldots,i_k,b\})\le c_1^{k+1} e^{-c_2 (b-a)}.
 \end{eqnarray}
Moreover, let $S^1,S^2$ be two independent copies of the copolymer, distributed according to the product measure
$\bP^{\otimes 2}_{N,\bgo}$. Then,
\begin{eqnarray}
  \label{eq:incontri}
  \bbE \, \bP^{\otimes 2}_{N,\bgo}\left(\nexists\,j:\, a <j< b, S^1_j=S^2_j=0\right)\le c_1
e^{-c_2(b-a)}.
\end{eqnarray}
\end{lemma}

\smallskip
\begin{remark}
\rm
It will be clear from the proof that the constant $c_2$
in \eqref{eq:ritorni} may be chosen smaller, but arbitrarily close to
$\mu(\bv)$. The quantitative estimate on the constant $c_2$ in \eqref{eq:incontri} 
that one can extract from the proof given below is instead substantially worse
and certainly not optimal.
\end{remark}
\bigskip

\noindent
{\it Proof of Theorem \ref{th:correlazioni}}.  
Let $\mathcal
S(A)=\{a_1,\ldots,a_2\}$ and $\mathcal S(B)=\{b_1,\ldots,b_2\}$. We
assume that $d(\mathcal S(A),\mathcal S(B))$ $>0$, otherwise
\eqref{eq:decad1} holds trivially with $c_1=2$. Without loss of
generality, we take $b_1>a_2$. Then, letting $\rm E$ be the event
$$
{\rm E}=\{\nexists\,j:\, a_2 <j< b_1,S^1_j=S^2_j=0\},
$$
one can write
\begin{eqnarray}
  \bE_{N,\bgo}(A B)-  \bE_{N,\bgo}(A)  \bE_{N,\bgo}(B)=
\bE^{\otimes 2}_{N,\bgo}\left\{\left(A(S_1)B(S_1)-A(S_1)B(S_2)\right)\ind_{\rm E}\right\}
\end{eqnarray}
since, by a simple symmetry argument based on the renewal property of $S$,
one can show that the above average vanishes, if conditioned to the complementary
of the event $\rm E$.
At this point, using  \eqref{eq:incontri} one obtains
\begin{eqnarray}
  \bbE\left[\left| \bE_{N,\bgo}(A B)-  \bE_{N,\bgo}(A)  \bE_{N,\bgo}(B)\right|\right]
\le 2 c_1 \Vert A\Vert_\infty\Vert B\Vert_\infty e^{-c_2(b_1-a_2)},
\end{eqnarray}
which is statement  \eqref{eq:decad1} of the theorem.

As for  \eqref{eq:bordo} we observe that, since we are assuming that $\mathcal S(A)\subset\{1,\ldots,k\}$,
one has the identity
\begin{equation}
\bE_{k, \bgo} \left(A \right)\, =\, 
\frac
{
\bE_{N, \bgo} \left(A\, \delta_k\right)
}
{
\bE_{N, \bgo} \left(\delta_k \right)
},
\end{equation}
where we recall that $\delta_k=\ind_{\{S_k=0\}}$.
Therefore, there exist
positive constants $c$ and $c^\prime$ such that
\begin{multline}
\bbE\left[
\left \vert
\bE_{N, \bgo}\left( A\right)-\bE_{k, \bgo}\left( A\right)
\right\vert
\right]
\, =\, 
\bbE\left[
\left \vert
\frac{
\bE_{N, \bgo}\left( A\right)\bE_{N, \bgo}\left( \delta_k\right)-
\bE_{N, \bgo}\left( A\,\delta_k\right)
}
{\bE_{N, \bgo}\left( \delta_k\right)
}
\right\vert
\right]
\\
 \le \,  c k^c
\bbE \left[\zeta(\gto_k)
\left\vert \bE_{N, \bgo}\left( A\right)\bE_{N, \bgo}\left( \delta_k\right)-
\bE_{N, \bgo}\left( A\,\delta_k\right)\right\vert 
\right]
\\
\le \, 
c k^c
 \left(\bbE \left[
\left\vert \bE_{N, \bgo}\left( A\right)\bE_{N, \bgo}\left( \delta_k\right)-
\bE_{N, \bgo}\left( A\,\delta_k\right)\right\vert
^2
\right]
\bbE\left[\zeta(\gto_k)^2
\right]\right)^{1/2}
\\
\le c^\prime k^{c} \Vert A \Vert_\infty
e^{-c_2 d(\mathcal S(A),\{k\})},
\end{multline}
where in the first inequality we have applied Lemma~\ref{th:lbpk},
in the second the Cauchy--Schwarz inequality and in the third 
Theorem~\ref{th:correlazioni}, formula \eqref{eq:decad1}.
Since $A$ is a local observable,  $d(\mathcal S(A),\{k\})= k (1+o(1))$ as $k \to \infty$
and therefore the proof of  \eqref{eq:bordo} is complete.
\smallskip

Finally,  \eqref{eq:exlim} is a consequence of the decay of the influence of 
boundary conditions expressed by  \eqref{eq:bordo}.
Note in fact that   \eqref{eq:bordo} states that  $\left\{\bE_{n, \bgo}\left( A\right)\right\}_n$
is a Cauchy sequence  in $L^1(\bbP)$. Therefore it convergence in $L^1$
toward a limit random variable that we denote $\bE_{\infty, \bgo}\left( A\right)$.
Therefore \eqref{eq:bordo} holds if we set  $N=\infty$ and, by the Fubini--Tonelli Theorem,  this
clearly implies that 
\begin{equation} 
 \bbE\left[\sum_{k} \left| \bE_{\infty,\bgo}(A)-\bE_{k,\bgo}(A)
\right|\right]<\infty,
\end{equation}
so the series in the expectation is $\bbP(\dd\bgo)$--a.s.
convergent, which implies the almost
sure convergence of $\left\{\bE_{n, \bgo}\left( A\right)\right\}_n$, 
that is  \eqref{eq:exlim}.
\hfill $\stackrel{\text {Theorem \ref{th:correlazioni}}}{\Box}$

\medskip

\noindent
{\sl Proof of Lemma \ref{lemma:ritorni}, Equation \eqref{eq:ritorni}}.
It is immediate to realize that, letting $i_0=a$ and $i_{k+1}=b$,
\begin{multline}
\label{eq:split}
\bbE\,\bP_{N,\bgo}(\{a \le j\le  b : S_j=0\}=\{i_0,i_1,\ldots,i_{k+1}\})
\\
\le\, \prod_{\ell=0}^{k} \bbE\left[
\frac{K(i_{\ell+1}-i_\ell)\zeta\left(\gto_{i_{\ell+1}}\right)}
{2 Z^\bv_{i_{\ell+1}-i_\ell,\theta^{i_\ell}\bgo}}
\left(1+
e^{-2\gl \sum_{j=i_\ell+1}^{i_{\ell+1}}(\go_j+h)}\right)
\right],
\end{multline}
where $K(\cdot)$ was defined in Section \ref{subs:fe} and $ \theta$ is the left shift.
Indeed, it suffices to apply Lemma \ref{prop:spezz} with $m=k+1$ and 
$A_j=\{S_i\ne 0, i_j< i< i_{j+1}\}$.
Thanks to Part 1 of Theorem \ref{th:zu} and to the exponential integrability of $\gto_1$, one has
$$
-\lim_{N\to\infty}\frac1N\log \bbE\left[
\frac{\zeta(\gto_N)}
{Z^\bv_{N,\bgo}}\left(1+
e^{-2\gl \sum_{j=1}^N(\go_j+h)}\right)
\right]\,=\,\mu(\bv) =: \mu > 0,
$$
so that one obtains from \eqref{eq:split}
\begin{equation}
\begin{split}
    \bbE\,\bP_{N,\bgo}(\{a \le j\le  b : S_j=0\}=\{i_0,i_1,\ldots,i_{k+1}\})
    \,\le\, &\prod_{\ell=0}^{k}c_3 e^{-\mu(i_{\ell+1}-i_\ell)/2}\\
\, =\, &c_3^{k+1}e^{-\mu(b-a)/2}.
\end{split}
\end{equation}

\hfill $\stackrel{\text{Lemma\; \ref{lemma:ritorni}, Eq.  \eqref{eq:ritorni}}}{\Box}$

\medskip

\begin{figure}[h]
\begin{center}
\leavevmode
\epsfxsize =12 cm
\psfragscanon
\psfrag{a}[c][l]{$a$}
\psfrag{b}[c][l]{$b$}
\psfrag{t10}[c][l]{$\eta^{(2)}_0$}
\psfrag{t11}[c][l]{$\eta^{(2)}_1$}
\psfrag{t12}[c][l]{$\eta^{(2)}_2$}
\psfrag{t13}[c][l]{$\eta^{(2)}_3$}
\psfrag{t14}[c][l]{$\eta^{(2)}_4$}
\psfrag{t15}[c][l]{$\eta^{(2)}_5$}
\psfrag{t16}[c][l]{$\eta^{(2)}_6$}
\psfrag{t20}[c][l]{$\eta^{(1)}_0$}
\psfrag{t21}[c][l]{$\eta^{(1)}_1$}
\psfrag{t22}[c][l]{$\eta^{(1)}_2$}
\psfrag{t23}[c][l]{$\eta^{(1)}_3$}
\psfrag{t24}[c][l]{$\eta^{(1)}_4$}
\psfrag{t25}[c][l]{$\eta^{(1)}_5$}
\psfrag{t26}[c][l]{$\eta^{(1)}_6$}
\psfrag{c10}[c][l]{$\chi^{(2)}_0$}
\psfrag{c11}[c][l]{$\chi^{(2)}_1$}
\psfrag{c12}[c][l]{$\chi^{(2)}_2$}
\psfrag{c13}[c][l]{$\chi^{(2)}_3$}
\psfrag{c14}[c][l]{$\chi^{(2)}_4$}
\psfrag{c15}[c][l]{$\chi^{(2)}_5$}
\psfrag{c16}[c][l]{$\chi^{(2)}_6$}
\psfrag{c20}[c][l]{$\chi^{(1)}_0$}
\psfrag{c21}[c][l]{$\chi^{(1)}_1$}
\psfrag{c22}[c][l]{$\chi^{(1)}_2$}
\psfrag{c23}[c][l]{$\chi^{(1)}_3$}
\psfrag{c24}[c][l]{$\chi^{(1)}_4$}
\psfrag{c25}[c][l]{$\chi^{(1)}_5$}
\psfrag{c26}[c][l]{$\chi^{(1)}_6$}
\psfrag{z}[c][l]{ $n$}
\psfrag{s}[c][c]{ $\dots$}
\psfrag{S1}[c][l]{ $S^2$}
\psfrag{S2}[c][l]{$S^1$}
\psfrag{n}[c][l]{$n$}
\epsfbox{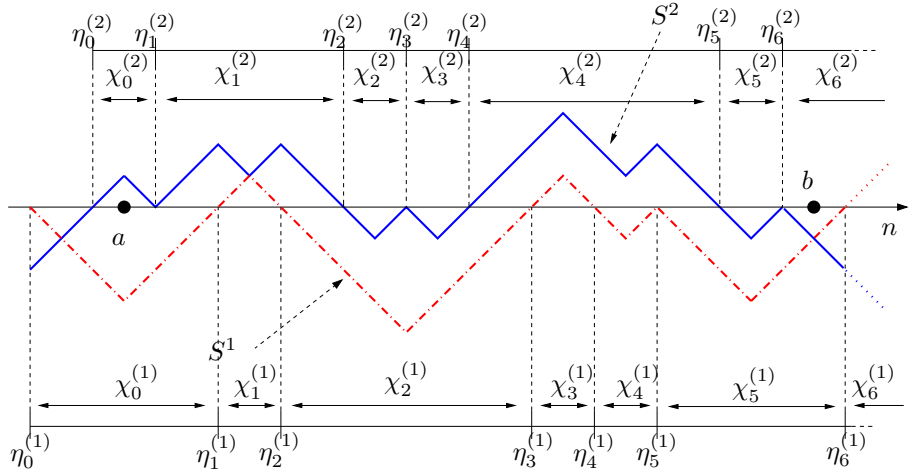}
\end{center}
\caption{\label{fig:figura} In the interval $\{a+1,\ldots,b-1\}$ the walk $S^1$ (dashed line) and the walk $S^2$ 
(full line) never 
 touch zero at the same time. We have marked the returns to zero $\eta^{(s)}_j$ and the 
lengths of the excursions $\chi^{(s)}_j$.
In this example, $r^{(1)}=5$ and $r^{(2)}=6$.
If for instance we choose $\kappa_2/\spC=4$, then in this example $S^1$ makes three {\sl short excursions},
and $S^2$ four.
Note that the two walks can meet and cross away from the line $S=0$.}
\end{figure}

\noindent
{\it Proof of Lemma \ref{lemma:ritorni}, Equation \eqref{eq:incontri}}.
In this proof the positive constants, typically dependent on $\bv$,
will be denoted by $d_0, d_1, \ldots$. The constants $c_1$ and $c_2$
are taken from \eqref{eq:ritorni}, but since $c_2$ is repeated several times
we set $\spC:=c_2$. 
Let us first define, for $s=1,2$,  $\eta^{(s)}_0=\sup\{0\le j\le a: S^s_j=0\}$, 
$\eta^{(s)}_i=\inf\{j> \eta^{(s)}_{i-1}:S^s_j=0\}$ for $i\ge 1$ and $r^{(s)}=\sup\{j: \eta^{(s)}_j< b\}$.
Then, for $j=0,\ldots,r^{(s)}$ we let $\chi^{(s)}_j=\eta^{(s)}_{j+1}-\eta^{(s)}_{j}\ge1$. 
We will refer to the interval $\{\eta^{(s)}_j,\cdots,\eta^{(s)}_{j+1}\}$ as to the $j^{th}$ excursion from zero 
of the walk $S^s$, and to 
$\chi^{(s)}_j$ as to its length, see Figure
\ref{fig:figura}. (Note that the $0^{th}$ and the $(r^{(s)})^{th}$ excursions may have an endpoint outside 
$\{a,\ldots,b\}$.)

The basic observation is that, as we will prove in a moment, there exists $\kappa_1,\kappa_2>0$ independent of $\bv$
such that
\begin{equation}
\label{eq:1/4}
\bbE\, \bP_{N,\bgo}\left(\sum_{0\le j\le r^{(1)}:\chi^{(1)}_j\ge \kappa_2 \spC^{-1}} 
\chi^{(1)}_j> \frac{b-a}4\right)\le 
d_1e^{-\kappa_1 \spC(b-a)},
\end{equation}
uniformly in $N\in \N$, $1\le a< b\le N$, for some finite constant $d_1:=d_1(\bv) >0$.
In words, this means that with high probability 
at least $3/4$ of the interval $\{a,\ldots,b\}$ is covered by excursions of $S$  whose length is smaller than
$\kappa_2\spC^{-1}$ (we will call them {\sl short excursions}).
To prove  \eqref{eq:1/4}, set $x=\chi^{(1)}_0$, $y=\chi^{(1)}_{r^{(1)}}$ and 
\begin{equation}
u:=\sum_{0<  j< r^{(1)}:\chi^{(1)}_j\ge \kappa_2\spC^{-1}}\chi^{(1)}_j,  
\end{equation}
so that the condition in the probability 
 in the left--hand side of \eqref{eq:1/4} reads 
$u+x \ind_{x\ge \kappa_2\spC^{-1}}+y\ind_{y\ge \kappa_2\spC^{-1}}\ge (b-a)/4$. 
Obviously, the number $\ell$ of excursions entirely contained in $\{a,\ldots,
b\}$ and 
of length at least $\kappa_2\spC^{-1}$ ({\sl long excursions}) is at most  
$\lfloor u\,\spC/\kappa_2\rfloor$, and one can (very roughly) bound above the
number of 
possible ways one can place them 
in the stretch $\{a,\ldots,b\}$ by
\begin{equation}
  \label{eq:combinatorio}
  \left(
    \begin{array}{c}
(b-a)\\
\ell
    \end{array}
\right)^2.
\end{equation}
These facts, together with a simple application of  Lemma \ref{prop:spezz} and  Eq. \eqref{eq:ritorni}, allow
to bound above the left--hand side of \eqref{eq:1/4} by
\begin{multline}
\label{eq:lunghe}
\sumtwo{x,y\ge 0,u\le (b-a):}{\left(x {\bf 1}_{x\ge \frac{\kappa_2}{\spC }}+y{\bf 1}_{y\ge \frac{\kappa_2}{\spC }}+
u\right)\ge \frac{b-a}4}  
(b-a)^2e^{-\spC  (u+x+y)} 
\sum_{\ell=0}^{\lfloor \frac{u\,\spC }{\kappa_2}\rfloor}\left(
    \begin{array}{c}
(b-a)\\
\ell
    \end{array}
\right)^2
c_1^{\ell+1}
\\
\le\,  d_ 0(b-a)^3 e^{- \spC  \frac{b-a}4}
\max(c_1, 1)^{\frac{\spC (b-a)}{\kappa_2}}\left(
    \begin{array}{c}
(b-a)\\
\lfloor\frac{\spC (b-a)}{\kappa_2} \rfloor
    \end{array}
\right)^2
\le d_1 e^{-\kappa_1\spC (b-a)},
\end{multline}
where the last inequality follows, if $\kappa_2$ is sufficiently large, from the Stirling formula. We stress that the constant $c_1$ is the one appearing in \eqref{eq:ritorni}.
The factor $(b-a)^2$ in \eqref{eq:lunghe} just takes care of the possible location of $\eta^{(1)}_1$ and 
$\eta^{(1)}_{r^{(1)}}$ in $\{a+1,\ldots,b-1\}$.

\begin{figure}[h]
\begin{center}
\leavevmode
\epsfxsize =12 cm
\psfragscanon
\psfrag{x1}[c][l]{$\bar x_1$}
\psfrag{x2}[c][l]{$\bar x_2$}
\psfrag{y1}[c][l]{$\bar y_1$}
\psfrag{y2}[c][l]{$\bar y_2$}
\psfrag{S1}[c][l]{$S^2$}
\psfrag{S2}[c][l]{$S^1$}
\psfrag{a}[c][l]{$a$}
\psfrag{p1}[c][l]{$p_1$}
\psfrag{p2}[c][l]{$p_2$}
\psfrag{p3}[c][l]{$p_3$}
\psfrag{p4}[c][l]{$p_4$}
\psfrag{p5}[c][l]{$p_5$}
\psfrag{b}[c][l]{$b$}
\psfrag{n}[c][l]{$n$}
\epsfbox{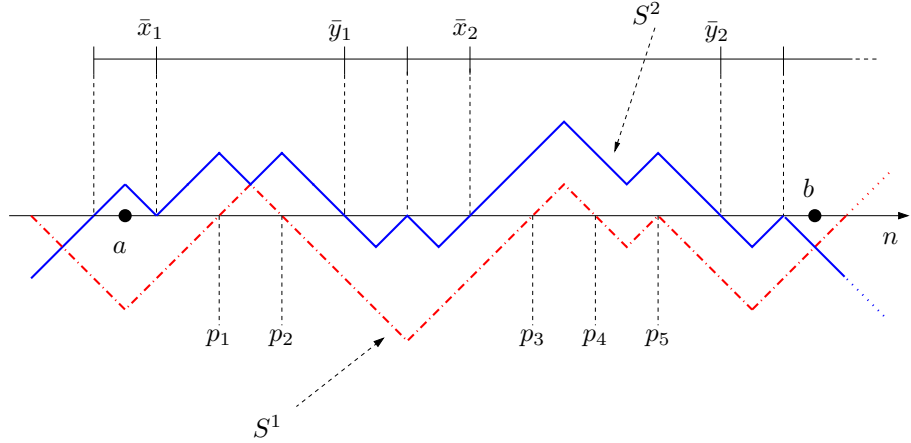}
\end{center}
\caption{\label{fig:trem} Consider again the example of Fig. \ref{fig:figura}.
If, say, $\kappa_2 \spC^{-1}=9$, then $\hat W=\{p_1,\ldots,p_5\}$. For $j=1,2$ one has
$(x_j,y_j)=(\bar x_1,\bar y_1)$, while $(x_j,y_j)=(\bar x_2,\bar y_2)$ if $j=3,4,5$.
}
\end{figure}

For ease of notation, let ${\rm B}_{a,b}$ be the event ${\rm B}_{a,b}=\{\nexists\,j:\, 
a <j< b,\, S^1_j=S^2_j=0\}$ that the two 
 walks do not touch zero at the same time 
 between $a$ and $b$, and for $s=1,2$ let ${\rm C}^{(s)}$ be the event
\begin{equation}
\label{eq:csm}
{\rm C}^{(s)}=
\left\{\sum_{0\le j\le r^{(s)}:\chi^{(s)}_j < \kappa_2\spC ^{-1}} \chi^{(s)}_j> \frac34(b-a)\right\}.  
\end{equation}
Then, from  Eq. \eqref{eq:1/4}, if 
$ \bar C={\rm C}^{(1)}\cap {\rm C}^{(2)}\cap{\rm B}_{a,b}$,
one has
\begin{equation}
\label{eq:intersezione}
\bbE \, \bP^{\otimes 2}_{N,\bgo}\left({\rm B}_{a,b}\right)\le 2 d_1
e^{-\kappa_1\spC (b-a)}+\,\bbE \, \bP^{\otimes 2}_{N,\bgo}\left(\bar C
\right).
\end{equation}
To estimate the last term in \eqref{eq:intersezione} let us notice that if the event $\bar C$ occurs then,
 denoting by $V^s$  the union of the short excursions of $S^s$, the set 
$V^1\cap V^2\cap\{a,\ldots,b\}$ contains at least $(b-a)/4$ sites. As a consequence, recalling that short excursions 
do not exceed $\kappa_2\spC^{-1}$ in length,  $ V:=\{j\in V^1\cap V^2\cap\{a,\ldots,b\}: 
S^1_j S^2_j=0\}$  contains at least $\lfloor \spC (b-a)/(8\kappa_2)\rfloor$ sites. 
In words, if $j\in V$ then either $S^1_j=0$ and $j$ belongs to a short excursion of $S^2$, or the
same holds interchanging the roles of $S^1$ and $S^2$. One can rewrite $V$ as the disjoint union
$V=W^1\cup  W^2$, where $W^s=\{j\in V: S^s_j=0\}$ and, of course, 
 at least one among $W^1$ and $W^2$ contains 
 $ \lfloor \spC(b-a)/(16\kappa_2)\rfloor$ points.
Therefore, using also the symmetry between $S^1$ and $S^2$, one has
\begin{equation}
\label{eq:proba}
\begin{split} 
 \bP^{\otimes 2}_{N,\bgo}\left(\bar C \right)
 \, &\le\,   2\, \bE^{\otimes 2}_{N,\bgo} \left( 
\ind_{\bar C\cap \{|W^1|\ge \lfloor \spC (b-a)/(16\kappa_2)\rfloor\}}\right)\\
 &\le \, 2\, \bE^{\otimes 2}_{N,\bgo} \left( 
\ind_{\bar C\cap\{|\hat W^1|\ge \lfloor \spC(b-a)/(16\kappa_2)\rfloor\}}\right),
\end{split}
\end{equation}
where $\hat W^1$ is the  subset of $\{a,\ldots,b\}$ 
which satisfies the following properties:
\smallskip
\begin{enumerate}
\item $S^1_j=0$ for every $j\in \hat W^1$;
\item for every $j\in \hat W^1$ there exist $a\le x_j<j<y_j\le b$ such  that 
  \begin{itemize}
  \item $S^2_{x_j}=S^2_{y_j}=0$
  \item $0<y_j-x_j<\kappa_2\spC^{-1}$
  \item $S^2_i\ne 0$ if $\{x_j<i<y_j$ and $S^1_i\ne0\}$.
  \end{itemize}
\end{enumerate}
\smallskip
Note that, since we are working on $\bar C$, 
$S^2_i\ne 0$ when $S^1_i=0$: this prescription is
not contained in the definition of  
$\hat W^1$.  
One can now write (see also Fig. \ref{fig:trem})
\begin{equation}
  \label{eq:condcomp}
  \begin{split}
\bP^{\otimes 2}_{N,\bgo}\left(\bar C\right)&\le
 2\, \bE_{N,\bgo}\left( \bE_{N,\bgo}^{\otimes 2}\left.\left(\ind_{\{\bar C\}} 
\ind_{\{|\hat W^1|\ge \lfloor \spC (b-a)/(16\kappa_2)\rfloor\}}\right|S^1
\right)\right)\\
&= 2\, \bE_{N,\bgo}\left( \sum_{\hat W}\bE_{N,\bgo}^{\otimes 2}
\left.\left(\ind_{\{\bar C\}} \ind_{\{\hat W^1=\hat W\}}
\right|S^1\right)\right)\\
&=  2\, \bE_{N,\bgo}\left(\ind_{\{\rm C^{(1)}\}} \sum_{\hat W}\bE_{N,\bgo}^{\otimes 2}
\left.\left(\ind_{\{\rm C^{(2)}\}}\ind_{\{\textrm{B}_{a,b}\}}
\ind_{\{\hat W^1=\hat W\}} \right|S^1\right)\right).
  \end{split}
\end{equation}
In the second step we have decomposed the probability by summing 
over all {\em a priori} admissible configurations $\hat W$ of the set $\hat W^1$, i.e., all possible subsets
of $\{a,\ldots,b\}$ containing at least $\lfloor \spC (b-a)/(16\kappa_2)\rfloor$ sites.
We are now going to relax 
the constraint given by 
$\ind_{\{\textrm{B}_{a,b}\}}$, but estimating the corresponding
ratio of probabilities. We claim in fact that, 
given $k,\ell\le \kappa_2 \spC^{-1}$ and $0<i_1<\ldots<i_k<\ell$,
we have
\begin{equation}
  \label{eq:cost}
  \frac{\bP_{\ell,\bgo}(  S_j\ne 0, 0<j<\ell)}
{\bP_{\ell,\bgo}( S_j\ne 0, j\notin\{0,i_1,\ldots,i_k, \ell\} )}\, \le\, 
 \prod_{r=1}^k\eta \left( \tilde \go _{i_r}
 \right),
\end{equation}
where
\begin{equation}
  \label{eq:def_eta}
  \eta(\gto)
\,  :=\, 
\left(1+d_2 \zeta\left(
\tilde  \go \right)
\right)^{-1},
\end{equation}
where $d_2$ is a positive constant that depends 
on $K(\cdot)$ and on the value of $\kappa_2 \spC^{-1}$.
The bound \eqref{eq:cost}, which will be applied with $S=S^2$,
is proven as follows.
First observe that
\begin{multline}
  \label{eq:cost1}
  \frac{\bP_{\ell,\bgo}(  S_j\ne 0, 0<j<\ell)}
{\bP_{\ell,\bgo}( S_j\ne 0, j\notin\{0,i_1,\ldots,i_k, \ell\} )}
\\
 \le\, 
 \frac{\bP_{\ell,\bgo}(  S_j\ne 0, 0<j<\ell)}
{\bP_{\ell,\bgo}( S_j > 0, j\notin\{0,i_1,\ldots,i_k, \ell\} )
+
\bP_{\ell,\bgo}( S_j < 0, j\notin\{0,i_1,\ldots,i_k, \ell\} )
},
\end{multline}
and in the ratio in the right--hand side 
we may factor the expression containing 
the copolymer energy term, that is we can set $\gl=0$
and we can restrict ourselves to considering $S\ge 0$.
We write:
\begin{equation}
\begin{split}
&\!\!\!\!\!\!\!\!\!\!\!\!\!\!\!\!\!\!\!\!\!\!  \, \frac
{
\bP_{\ell,\bgo}
\left( S_j> 0, j\notin\{0,i_1,\ldots,i_k, \ell\} \right)
}
{\bP_{\ell,\bgo}\left(  S_j > 0, 0<j<\ell\right)}
\\
&=\, 
\sum_{A \subset \{i_1, \ldots, i_k\}}
\frac
{
\bP_{\ell,\bgo}\left( 
S_j= 0, j\in A, \, 
S_j> 0, j\in \{1,\ldots, \ell -1\} 
\setminus A 
\right)
}
{K(\ell)/2}
\\
&\ge \, 
\sum_{A \subset \{i_1, \ldots, i_k\}} 
\prod_{j\in A} \epsilon_\ell  \zeta \left( \tilde \go _j \right)\, =\,
\prod_{j \in \{i_1, \ldots, i_k\}} 
\left( 1+\epsilon_\ell \zeta \left( \tilde \go _j \right) \right),
\end{split}
\end{equation}
where $\epsilon_\ell$ may be chosen, with a very rough
estimate, equal to $\min_{j =1, \ldots, \ell}K(j)^2/4$.
Since $\ell\le \kappa_2 \spC^{-1}$, we obtain
\eqref{eq:cost}.

Therefore from
\eqref{eq:condcomp}, using \eqref{eq:cost},
we can extract
\begin{equation}
  \label{eq:tremenda}
\bP^{\otimes 2}_{N,\bgo}\left(\bar C\right)\, \le \,
  2\, \bE_{N,\bgo}\left(\ind_{\{\rm C^{(1)}\}} \sum_{\hat W}\bE_{N,\bgo}^{\otimes 2}
  \left.\left(\ind_{\{\rm C^{(2)}\}}
\ind_{\{\hat W^1=\hat W\}} \right|S^1\right)\prod_{j\in \hat W}\eta(\gto_j)\right).
\end{equation}

The remaining problem now is that $\eta( \tilde \go_n)<1$ is random and 
can get arbitrarily close to $1$. This however does not happen too often: 
 in fact there exists $\kappa_3$  such that, for every $\bv\in \cL$, 
\begin{equation}
  \label{eq:probU}
\bbP\left(U\right)\, := \,  
\bbP\left(\left|\left\{a<i<b: \omega_i<-\kappa_3\right\}\right|
\ge\left\lfloor\frac{\spC(b-a)}{32\kappa_2}\right\rfloor
\right)\le d_3 e^{-d_4(b-a)},
\end{equation}
for every value of $(b-a)$.
The bound \eqref{eq:probU} follows from a direct (large deviation) estimate on
the binomial random variable with  parameters
$ p:=\bbP\left(\tilde \go_i<- \kappa_3 \right)$, which
can be made arbitrarily small, and $b-a$.
It suffices for example that $2p\le \spC/32 \kappa_2$.

Thanks to  $|W|\ge \lfloor \spC (b-a)/(16\kappa_2)\rfloor$,   \eqref{eq:tremenda}
 implies that, on the complementary of $U$,
\begin{equation}
  \label{eq:compU}
  \begin{split}
\bP^{\otimes 2}_{N,\bgo}\left(\bar C\right)\le (1+d_2 
\zeta(-\kappa_3))^{-\frac{\spC(b-a)}{32\kappa_2}}\le d_5 e^{- d_6
(b-a)}.
  \end{split}
\end{equation}
Together with  \eqref{eq:intersezione} and \eqref{eq:probU}, this completes the proof of
 \eqref{eq:incontri}.
\hfill $\stackrel{\text{Lemma\; \ref{lemma:ritorni},   \eqref{eq:incontri}}}{\Box}$

\section{Regularity of the free energy}
\label{sec:reg}

\noindent
{\it Proof of  Theorem \ref{th:cinf}}. 
Thanks to the Ascoli-Arzel\`a Theorem, it is sufficient to show that, for every integer $k$, the
$k^{th}$ derivative of $(1/N)\bbE\log Z^\bv_{N,\bgo}$ with respect to any of the parameters $\gl,h,\tilde\gl,\tilde h$
is bounded above uniformly in $N$.
For definiteness, let us show this property for 
\begin{eqnarray}
  \frac{\partial^k}{\partial{\tilde\gl^k}} \frac1N \bbE \log Z^{\bv}_{N,\bgo}.
\end{eqnarray}
The above derivative is given by 
\begin{eqnarray}
\label{eq:derivk}
 \frac{\partial^k}{\partial{\tilde\gl^k}} \frac1N \bbE \log Z^{\bv}_{N,\bgo} =
 \frac1N \sum_{1\le n_1,n_2,\ldots,n_k\le N}\bbE\left\{ \gto_{n_1}\ldots \gto_{n_k}
\bE_{N,\bgo}(\delta_{n_1};\delta_{n_2};\ldots;\delta_{n_k})\right\},
\end{eqnarray}
which is expressed through truncated correlation functions (Ursell functions) defined as 
\begin{eqnarray}
  \label{eq:ursell}
 \bE_{N,\bgo} \left(A_1;\ldots;A_k\right)=
\sum_{\mathcal P}(-1)^{|\mathcal P|-1}(|{\mathcal P}|-1)!\prod_{P\in {\mathcal P}}
 \bE_{N,\bgo}\left( \prod_{p\in P}A_p\right),
\end{eqnarray}
where the sum runs over all partitions ${\mathcal P}$ of $\{1,\ldots,k\}$ into subsets $P$.
Starting from the property \eqref{eq:decad1} of decay of correlations for every
pair of bounded local observables, one can prove by induction over $k\ge2$ that 
\begin{eqnarray}
\label{eq:ursell2}
\!\!\!\!\!\!\!\!\!\!\!\!\!\! &&
\bbE \left[\left|\bE_{N,\bgo}(A_1;\ldots;A_k)\right|\right] \\
\!\!\!\!\!\!\!\!\!\!\!\!\!\! && \qquad
\le c_1^{(k)}\Vert A_1\Vert_\infty\ldots \Vert A_k\Vert_\infty
\exp(-c_2^{(k)}d(\mathcal S(A_1),\ldots,\mathcal S(A_k))),
\nonumber
\end{eqnarray}
for some finite 
positive constants $c_1^{(k)},c_2^{(k)}$ where, if $I$ is the smallest interval including the supports
of $A_1,\ldots,A_k$, then $d(\mathcal S(A_1),\ldots,\mathcal S(A_k))
\equiv |I|-\sum_{\ell=1}^k|\mathcal S(A_\ell)|$.
A proof that the exponential decay of average two--points correlations implies exponential decay of the average 
$n$--point 
truncated correlations can be found
for instance in \cite{cf:griffiths},
in the context of disordered $d$--dimensional spin systems in the high temperature or 
large magnetic field regime (see Remark~\ref{rem:g} below for a sketch
of the proof). In \cite{cf:griffiths}, explicit bounds for the constants
$c_1^{(k)},c_2^{(k)}$ are also given, which is not needed in our case. The proof in \cite{cf:griffiths}, which is 
basically an application of H\"older's inequality,  can be transposed
to the present context almost without changes, to yield \eqref{eq:ursell2}.
Then, after an application of the 
Cauchy-Schwarz inequality, 
it is immediate to realize that the sum in \eqref{eq:derivk} converges uniformly in $N$, for every $k$.

\hfill $\stackrel{\text{Theorem\;\ref{th:cinf}}}{\Box}$

\medskip

\begin{remark}\label{rem:g}\rm
A sketch  of how  \eqref{eq:ursell2} is deduced from \eqref{eq:decad1}
goes as follows. 
Consider the case $k=3$ and assume that 
$A_i=\delta_{n_i}$ (which is just the case we need in view of  \eqref{eq:derivk}) and, without loss of generality,
let $n_1\le n_2\le n_3$. Then consider the simple identities
\begin{equation}
\begin{split}
\bE_{N,\bgo}\big(&\delta_{n_1};\delta_{n_2};\delta_{n_3}\big)
\\
&=\,  
\bE_{N,\bgo}(\delta_{n_1}\delta_{n_2};\delta_{n_3})-\bE_{N,\bgo}(\delta_{n_1}) \bE_{N,\bgo}(\delta_{n_2};\delta_{n_3})-
\bE_{N,\bgo}(\delta_{n_2}) \bE_{N,\bgo}(\delta_{n_1};\delta_{n_3})\\
&=\, \bE_{N,\bgo}(\delta_{n_2}\delta_{n_3};\delta_{n_1})-
\bE_{N,\bgo}(\delta_{n_2}) \bE_{N,\bgo}(\delta_{n_1};\delta_{n_3})-
\bE_{N,\bgo}(\delta_{n_3}) \bE_{N,\bgo}(\delta_{n_1};\delta_{n_2}).
\end{split}
\end{equation}
From the first identity and  \eqref{eq:decad1} we obtain 
\begin{eqnarray}
  \bbE \left[\left|\bE_{N,\bgo}(\delta_{n_1};\delta_{n_2};\delta_{n_3})\right|\right]
\le 3c_1 e^{-c_2 (n_3-n_2)},
\end{eqnarray}
while from the second identity we obtain the bound 
$3c_1 \exp(-c_2 (n_2-n_1))$ on the same quantity. Therefore 
\begin{eqnarray}
   \bbE \left[\left|\bE_{N,\bgo}(\delta_{n_1};\delta_{n_2};\delta_{n_3})\right|\right]\le 
3c_1e^{-c_2((n_3-n_2)+(n_2-n_1))/2},
\end{eqnarray}
which is just the statement of \eqref{eq:ursell2} in this specific case.
\end{remark}

\section{On the maximal excursion}
\label{sex:exc}

\subsection{Proof of Theorem \ref{th:zu}}
\label{sec:maximal}
\noindent 
{\it Part 1.}
The existence of the limit follows from the subadditivity of
$\left \{ \log \bbE \left[ (1+\exp(-2\lambda\sum_{n=1}^N(\go_n+h))/ Z_{N, \bgo}^{\bv}
\right]\right\}_N$, which is an immediate consequence
of the renewal property (of $\bP$) and of the IID property (of $\bgo$):
for every $M \in \N$, $M<N$, we have in fact
\begin{equation}
\label{eq:subaddmu}
\begin{split}
\bbE \left[  \frac{1+e^{-2\lambda\sum_{n=1}^N(\go_n+h)}}
{ Z_{N, \bgo}^\bv}
\right]
\, &\le \, 
\bbE \left[ \frac{(1+e^{-2\lambda\sum_{n=1}^M(\go_n+h)})(1+e^{-2\lambda\sum_{n=M+1}^N(\go_n+h)})}
 {Z_{M, \bgo}^\bv Z_{N-M, \theta^M \bgo}^\bv} 
\right]\\
&=\bbE \left[  \frac{1+e^{-2\lambda\sum_{n=1}^M(\go_n+h)}}
{ Z_{M, \bgo}^\bv}
\right]
\bbE \left[  \frac{1+e^{-2\lambda\sum_{n=1}^{N-M}(\go_n+h)}}
{ Z_{N-M, \bgo}^\bv}
\right].
\end{split}
\end{equation}
The inequality $\mu (\bv)\ge 0$ is an immediate
consequence of the elementary lower bound
\begin{equation}
Z_{N, \bgo} \ge \zeta(\gto_N) \left(1+e^{-2\lambda\sum_{n=1}^N(\go_n+h)}\right)
\bP\left( \gO_N^+\right).
\end{equation}
A more refined lower bound on $\mu (\bv)$, valid in the localized region, 
follows from the concentration inequality: 
call $E_N$ the event 
$$
E_N=\left\{ \bgo:\,-(4\gl/N)\sum_{n=1}^N(\go_n+h) 
< \tf (\bv )/2 <  (1/N)\log Z_{N, \bgo} ^{\bv} \right\}.
$$
Since for $N$ sufficiently large 
$\bbP \left( E^\complement_N\right) \le \kappa_1\exp (-\kappa_2 N \tf(\bv)^2
/\max(\gl, \tilde \gl)^2)$ with $\kappa_1$ and $\kappa_2$ suitable positive constants, 
one has
\begin{eqnarray*}
\!\!\!\!\!\!\!\!\!\!\!\!\!\!\! &&
\bbE \left[ \frac{1+e^{-2\gl\sum_{n=1}^N(\go_n+h)}} 
{Z_{N,\bgo}^\bv}\right] \, \le \, 
2\exp\left(-N\tf (\bv)/4\right)+ 
\frac
{\bbE \left[\zeta(\gto_N)^{-1}\ind_{\{ \, E_N^\complement\}}
\right]}{\bP\left( \gO_N^+\right)} \\
\!\!\!\!\!\!\!\!\!\!\!\!\!\!\! && \qquad\qquad\qquad\qquad\qquad
\le 2\exp\left(-N\tf (\bv)/4\right)+ \kappa'_1\exp
\left(-\kappa'_2 N\tf(\bv)^2 /\max(\gl, \tilde \gl)^2\right),
\end{eqnarray*}
which immediately implies $\mu(\bv)>0$ and, for $\tf(\bv)$ sufficiently small,
$2\mu(\bv)>\kappa'_2\tf(\bv)^2/\max(\gl, \tilde \gl)^2$.

Finally, $\mu(\bv) \le \tf (\bv)$ is an immediate consequence
of Jensen's inequality.

\hfill $\stackrel{\text{Theorem\;\ref{th:zu}\; part 1}}{\Box}$

\bigskip

\noindent
{\it Part 2.}
Throughout this proof $\bv \in \cL$, so that $\mu(\bv)>0$, and we set
 $a_N^\pm := \lfloor (1\pm \gep) \log N/\mu( \bv) \rfloor$.
We start with the proof of
\begin{equation}
\label{eq:ubexc}
\lim_{N\to \infty}
\bbE\, \bP_{N, \bgo}\left(  {\boldsymbol \Delta_N} > a_N^+\right) \, =\, 0,
\end{equation}
which clearly implies \eqref{eq:problb}.
For $n=0,1,2, \ldots$ we set
\begin{equation}
\label{eq:En}
E_n\, :=\, \left\{S: \, S_n=0, \, S_{n+j}\neq 0 \text{ for } j=1,2, \ldots, a_N^+\right\},
\end{equation}
so that $\{\boldsymbol \Delta_N > a_N^+ \} = \cup_n E_n$, $n $ ranging up to $N-a_N^+-1$,
and we have
\begin{equation}
\label{eq:unionb}
\bbE\, \bP_{N, \bgo}\left( {\boldsymbol \Delta_N} > a_N^+\right)\, \le \,
\sum_{n:\,  n \le N-a_N^+-1}
\bbE\, \bP_{N, \bgo}\left( E_n \right) .
\end{equation}
Let us estimate the terms in the sum by writing first $E_n$
as the disjoint union of the sets
\begin{equation}
\label{eq:Enl}
E_{n, \ell}\, :=\, \left\{S: \, S_n=0, \, S_{n+j}\neq 0 \text{ for } j=1,2, \ldots, \ell -1, \,
S_{n+\ell }=0\right\},
\end{equation}
with $\ell \in \N$ and $ a_N^+<\ell\le N-n$.
Recalling the notations of Section \ref{sec:conv}, we have the bound 
\begin{equation}
\label{eq:bdr1}
\begin{split}
\bP_{N, \bgo}\left( E_{n,\ell}\right)\, 
& \le \, 
 \frac{Z^{\bv}_{N, \bgo}\left(E_{n,\ell}\right)}
{Z^{\bv}_{N, \bgo}\left(S_n=S_{n+\ell}=0\right)}
\\
&\le\, 
\frac{1+ \exp\left(-2\gl \sum_{j=1}^\ell \left( (\theta^n\go)_j +h\right)\right)
} 
{
2
Z^{\bv}_{\ell, \theta^n \bgo}
} \zeta(\gto_{n+\ell}).
\end{split}
\end{equation}
This is an immediate consequence of the renewal property of $S$.
Notice that once we take the expectation of both sides
of \eqref{eq:bdr1} we may set $n=0$ in the right--hand side
and therefore,
by \eqref{eq:mu}, for every $n$ and every $\ell\ge \ell _0$, $\ell_0$ sufficiently large,
we have that 
\begin{equation}
\label{eq:expo}
\bbE\, \bP_{N, \bgo}\left( E_{n,\ell} \right)
\le
\exp\left(-\ell \mu (\bv ) (1-(\gep/2))
\right).
\end{equation}
Indeed, the factor $\zeta(\gto_{n+\ell})=\exp (\tilde \gl(\tilde \go_{n+\ell} + \tilde h))$ 
is negligible for $\ell$ large since,
with  probability at least of order  $1-O(\exp(-c\ell^{3/2})$,  $|\gto_{n+\ell}|$  does not exceed 
$\ell^{3/4}$, see  \eqref{eq:Lip}, so that it does not modify the exponential behavior \eqref{eq:expo}.

Therefore, for $N$ sufficiently large and $\gep$ sufficiently small
\begin{equation}
\begin{split}
\bbE\, \bP_{N,\bgo} \left(E_n\right)\, &=\,
\sumtwo{\ell \in \N: }{ a_N^+<\ell\le N-n}
\bbE\,\bP_{N, \bgo}\left( E_{n,\ell} \right)
\\ 
&\le\, 
\sum_{\ell >a^+_N} \exp\left(-\ell  \mu (\bv ) (1-(\gep /2))\right)  \le 
N^{-1- (\gep /4)}.
\end{split}
\end{equation}
Going back to
\eqref{eq:unionb} we see that $\bbE\, \bP_{N, \bgo}\left( {\boldsymbol \Delta_N} > a_N^+\right)\le N^{-\gep/4}$ and 
\eqref{eq:ubexc} is proven.

\medskip
\begin{remark}
\label{rem:cutoff}
\rm
It is immediate to realize by looking at the proof that
\begin{equation}
\bbE\,\bP_{N, \bgo}\left(
\boldsymbol \Delta _N \, >\, 
C \log N \right)\le N^{-L(C)},
\end{equation}
with $L(C)\nearrow \infty$ as $ C\nearrow \infty$.
\end{remark}

\medskip

Let us then turn to proving 
\begin{equation}
\label{eq:lbexc}
\lim_{N\to \infty}
 \bP_{N, \bgo}\left(  {\boldsymbol \Delta_N} \ge  a^-_N\right) \, =\, 1, \ \ \ \  \bbP (\dd \bgo)-\text{a.s.}.
\end{equation}
Let us set 
$j_N:= \lfloor N /(\log N)^2 \rfloor -1$ and $n_j = \left\lfloor Nj/(j_N+1)\right\rfloor$ for $j=1, 2 , \ldots, j_N$.

\begin{figure}[h]
\begin{center}
\leavevmode
\epsfxsize =11.5 cm
\psfragscanon
\psfrag{0}[c][l]{$0$}
\psfrag{N}[c][l]{$N$}
\psfrag{n1}[c][l]{$n_1$}
\psfrag{n2}[c][l]{$n_2$}
\psfrag{n3}[c][l]{$n_3$}
\psfrag{n4}[c][l]{$n_4$}
\psfrag{n5}[c][l]{$n_5$}
\psfrag{an}[c][l]{$a_N^-$}
\epsfbox{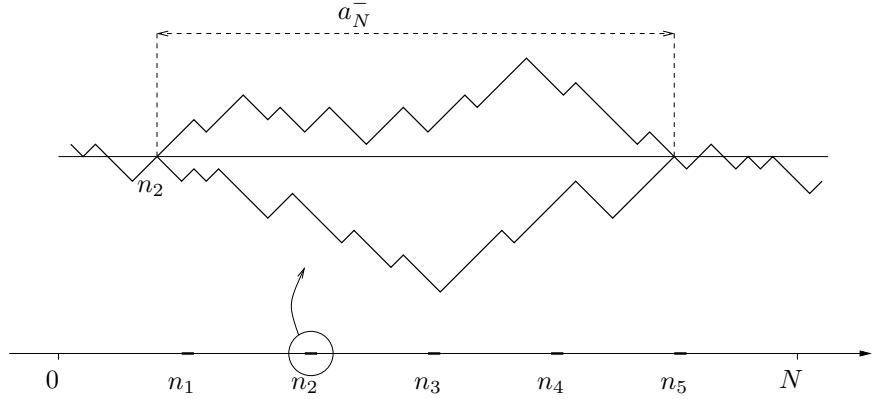}
\end{center}
\caption{\label{fig:LB} 
The lower bound on the length of the maximal excursion is achieved by focusing 
on what happens right after $j_N= N(1+o(1))/(\log N)^2$ 
sites $n_j$,  equal to $n_j=j (\log N)^2$ (modulo lattice corrections). In the Figure, $j_N=5$.
One looks at the probability of finding the (positive or negative) excursion that starts
exactly at $n_j$ and comes back to zero at $n_j+a_N^-$,
 at least for one $j$. Note that $a_N^-= (1-\gep)(\mu(\bv))^{-1}(1+o(1))\log N$ is much smaller
 than $n_j-n_{j-1}$ and this allows decoupling of these occurrences. }
\end{figure}

We now consider the family of events $\left\{E_{n_j, a_N^-}
\right\}_{j=1, 2 , \ldots, j_N}$, recall the definition of $E_{n,\ell}$ in
\eqref{eq:Enl}, and we observe that $\cup_j E_{n_j, a_N^-} \subset \{ \boldsymbol \Delta_N \ge a_N^-\}$.
In words,  we are simply going to compute the probability that
the walk makes at least one  excursion of length exactly equal to $a_N^-$ at
$j_N$ prescribed locations: the excursions have to start at $n_j$ for some $j$, see Figure \ref{fig:LB}. 
Therefore
\begin{equation}
\label{eq:RNdef}
 \bP_{N, \bgo}\left(  {\boldsymbol \Delta_N} <  a^-_N\right)
 \le \bP_{N, \bgo}\left( \cap_{j=1}^{j_N} E_{n_j, a_N^-}^\complement \right)
 =
 \prod_{j=1}^{j_N} \left(1- \bP_{N, \bgo}\left(E_{n_j, a_N^-} \right)\right)+ R_N(\bgo)
 ,
\end{equation}
where the last step defines $R_N(\bgo)$ (we anticipate that, by Theorem~\ref{th:correlazioni}, formula \eqref{eq:decad1},
  $R_N(\bgo)$
is negligible,
details are postponed to 
Lemma~\ref{th:RN} below).
We need therefore a lower bound on
$\bP_{N,\go}\left(E_{n_j, a_N^-} \right)$: we will find a lower bound
on this quantity that depends on $\bgo$ only via $\go_n$ and $\tilde \go_n$
with $n=n_j,n_j +1, \ldots, n_j+a_N^-$ and this will allow the direct use of
independence when taking the expectation with respect to $\bbP$. 
We use the explicit formula
\begin{multline}
\label{eq:rtio}
\bP_{N, \bgo}\left(E_{n_j, a_N^-} \right)\, =
\\
\frac
{Z^{\bv}_{n_j,\bgo} K\left(a_N^-\right) 
\zeta(\gto_{n_j+a^-_N})
\left(1+ \exp\left(-2\gl \sum_{n=n_j+1}^{n_j+a_N^-}(\go_n+h)\right)\right)
Z^{\bv}_{N-n_j-a_N^-, \theta^{n_j+a_N^-}\bgo}
}
{
2Z^{\bv}_{N, \bgo}
}.
\end{multline}

We proceed by finding an upper bound on the denominator
and we first observe that,
by \eqref{eq:ubexc},
the event $A_{N}:=\left\{ \bgo:\,
 Z^{\bv}_{N, \bgo}(\boldsymbol  \Delta_N \le C \log N)/
Z^{\bv}_{N, \bgo}\ge 1-\gd\right\}$ has $\bbP(\dd \bgo)$--probability
tending to $1$ for $C>1/\mu (\bv )$ and $\gd>0$:
this is simply due to the fact that 
$Z^{\bv}_{N, \bgo}(\boldsymbol \Delta_N \le C \log N)/
Z^{\bv}_{N, \bgo}= \bP_{N, \bgo}( \boldsymbol\Delta_N  \le C \log N)$.
 We will actually choose $C$ larger,
so that to guarantee that $\bbP[A_N] \ge 1-N^{-2}$ for $N $ large, cf. Remark~\ref{rem:cutoff}.
Of the requirements defining the event $\boldsymbol  \Delta_N \le C \log N$
we now take advantage
only of the fact that there exists a return to zero at distance at most $C\log N$ from both sites
$n_j$ and $n_j+a^-_N$,
for every $j$.
By Lemma~\ref{th:lbpk} we have that there exists $c>0$ such that for $\bgo \in A_{N}$
\begin{equation}
\begin{split}
Z^{\bv}_{N, \bgo}\, & \le\, \frac 1{1-\gd} Z^{\bv}_{N, \bgo}\left(\min_i |\tau_i-n_j|, \min_i |\tau_i-(n_j+a^-_N)|
\le C \log N\right)
\\
 & \le\,  
c (\log N )^c Z^{\bv}_{n_j,\bgo}Z^{\bv}_{a_N^-, \theta^{n_j}\bgo}
Z^{\bv}_{N-n_j-a_N^-, \theta^{n_j+a_N^-,}\bgo}\zeta(|\gto_{n_j}|)\zeta(|\gto_{n_j+a^-_N}|),
\end{split}
\end{equation}
   so that from \eqref{eq:rtio} there exist $c ,c^\prime 
>0$  such that
\begin{equation}      
   \label{eq:rtio2}
\bP_{N, \bgo}\left(E_{n_j, a_N^-} \right)\,\ge \,
c(\log N)^{-c^\prime} \frac{
1+ \exp\left(-2\gl \sum_{n=n_j+1}^{n_j+a_N^-}(\go_n+h)\right)}
{\zeta(|\gto_{n_j}|)\zeta(|\gto_{n_j+a^-_N}|)\,Z^{\bv}_{a_N^-, \theta^{n_j}\bgo}
}
=: \, Q_{N, j}(\bgo)
,
\end{equation}
where the inequality holds 
for $N $ sufficiently large. Notice that the
 random variables $\left\{Q_{N, j}(\bgo) \right\}_j$
are independent,
since $Q_{N, j}(\bgo) $ depends
 on $(\go_n, \tilde\go_n)$ 
 with $n$ only in $\{n_j, n_j +1, \ldots, n_j +a_N^-\}$.
Going back to \eqref{eq:RNdef}, we have
\begin{equation}
\label{eq:nowE}
\begin{split}
\bP_{N, \bgo}\left(  {\boldsymbol \Delta_N} <  a^-_N\right)\, &=\,\bP_{N, \bgo}\left(  {\boldsymbol \Delta_N} <  a^-_N\right)\left(\ind_{A_N}
(\bgo)+
\ind_{A_{N}^\complement}(\bgo)\right)\\
&\le \, 
\prod_{j=1}^{j_N}
\left( 1- Q_{N, j}(\bgo) \right)\, +\,  \ind_{A_{N}^\complement}(\bgo) \, +\, R_N(\bgo).
\end{split}
\end{equation} 
We take now  the $\bbP(\dd \bgo)$--expectation of
both sides of \eqref{eq:nowE}: by independence, by the choice of $C$ in the definition
of $A_N$ and by Lemma~\ref{th:RN}
\begin{equation}
\label{eq:dgfd}
\begin{split}
\bbE\, \bP_{N, \bgo}\left(  {\boldsymbol \Delta_N} <  a^-_N\right)\, 
&\le\,
\prod_{j=1}^{j_N}
\left( 1- \bbE \left[Q_{N, j}(\bgo) \right]\right)
\, +\,  
\bbP\left[A_{N}^\complement\right] \, +\, \bbE
\left[\left\vert R_N(\bgo)\right\vert\right]
\\
&=\,
\prod_{j=1}^{j_N}
\left( 1- \bbE \left[Q_{N, j}(\bgo) \right]\right)\, +\, O(1/N^2).
\end{split}
\end{equation}
We are reduced to estimating
$\bbE[Q_{N, j}(\bgo) ]$.
 By \eqref{eq:mu}, for $\gep$ sufficiently small and $N$ sufficiently large, we have 
\begin{equation}
\bbE\left[Q_{N, j} (\bgo)\right]\,\ge \, 
\exp\left(- \mu (\bv ) (1+\gep/2) a_N^-\right)\,\ge\, N^{-1+ \gep/4}. 
\end{equation}
From this  one gets
\begin{equation}
\label{eq:cfrn}
 \bbE\,\bP_{N, \bgo}\left(  {\boldsymbol \Delta_N} <  a^-_N\right)
 \,\le\, 
 \left( 1-N^{-1+ \gep/4}\right)^{j_N} \, +\,
O(1/N^2)\, = \, O(1/N^2).
\end{equation}
By applying the Markov inequality and Borel--Cantelli lemma
we conclude the proof of \eqref{eq:lbexc}.

\hfill $\stackrel{\text{Theorem\;\ref{th:zu}\; part 2}}{\Box}$

\medskip

\begin{lemma}
\label{th:RN}
For  every $m>0$ we have
\begin{equation}
  \lim_{N\to \infty}N^m\,\bbE\left[\left\vert R_N(\bgo) \right\vert\right] \, =\, 0,
\end{equation}
with $R_N(\bgo)$ defined in \eqref{eq:RNdef}.
\end{lemma}
\medskip

\noindent
{\it Proof.}
Observe that 
a direct application of Theorem~\eqref{th:correlazioni}, formula
\eqref{eq:decad1}, yields that for $k=2, 3, \ldots, j_N$ we have
\begin{multline}
\bbE\left[
\left\vert
\bP_{N,\bgo} \left( \cap_{j=1}^k E^\complement_{n_j,a^-_N}\right)-
\bP_{N,\bgo} \left( \cap_{j=1}^{k-1} E^\complement_{n_j,a^-_N}\right)\bP_{N,\bgo} \left( E^\complement_{n_k,a^-_N}\right)
\right\vert
\right]
\\
\le\,  c_1 \exp\left(-c_2 \frac{(\log N)^2}2\right)
.
\end{multline}
One now applies iteratively this inequality starting from $k=j_N$, down to
$k=2$, obtaining that $\bbE\left[\left\vert R_N(\bgo) \right\vert\right]$
is bounded above by $
c_1 N \exp\left(-c_2 (\log N)^2/2\right)$.

\hfill $\stackrel{\text{Lemma\; \ref{th:RN}}}{\Box}$

\subsection{Proof of Proposition \ref{th:as}}
It makes use of the identity
\begin{multline}
  \label{eq:id1}
  \bP_{N, \bgo}\left( \boldsymbol \Delta_N (k) = s\right) = 
  \\
  K(s)
\sumtwo{l\in \N\cup \{0\},r\in \N :}{l+r=s, \, k-l\ge 0, \, k+r\le N}
\!\!\!\!\!
\frac{Z^{\bv}_{k-l, \bgo} 
\left(
1+ e^{ -2 \gl \sum_{n=k-l+1}^{k+r} (\go_n+h)}\right)\zeta(\gto_{k+r})Z^{\bv}_{N-k-r,\theta^{k+r}\bgo}}
{2Z^{\bv}_{N,\bgo}}.
\end{multline}
\smallskip

\noindent
{\it Upper bound.}
We first observe that, reasoning as  in Section \ref{sec:maximal}, we have
\begin{equation}
\label{eq:s1as}
\bP_{N, \bgo}
\left( \boldsymbol \Delta_N (k) = s\right) \, \le \,
\sumtwo{l\in \N\cup \{0\},r\in \N :}{l+r=s}
G_{l,r}\left(\theta^k \bgo \right),
\end{equation}
where
\begin{equation}
G_{l,r}\left(\bgo \right)\, =\,
K(l+r)
\frac{
1+ \exp \left( -2 \gl \sum_{n=-l+1}^r (\go_n+h)\right)
}
{2Z_{l+r, \theta^{-l}\bgo}
}\zeta(\gto_{r})
,
\end{equation}
(recall that we are working here with two--sided disorder sequences, so that the $\bgo$ variables
may have negative indices).
Now we claim that for every $\gep>0$ there exists $s_0:= s_0\left(\bgo\right)$, $\bbP
(\dd \bgo)$--a.s.
finite, such that
\begin{equation}
\label{eq:asclaim}
G_{l,r}\left(\bgo \right)\, \le \, \exp \left( - s (\tf (\bv)-\gep/2)\right),
\end{equation}
for every $l$ and $r$ such that $l+r = s$ and every $s\ge s_0$.
Of course \eqref{eq:asclaim}
implies that for every $\gep>0$
\begin{equation}
\label{eq:asC}
C^{(2)}_\gep (\bgo)\, :=\,
\sum_{s \in  \N }
\sumtwo{l\in \N\cup \{0\},r\in \N :}{l+r=s}
G_{l,r}\left(\bgo \right) \exp \left(  s (\tf (\bv)-\gep)\right) \, < \infty, \ \ \
\bbP(\dd \bgo)-\text{a.s.}.
\end{equation}
By combining \eqref{eq:s1as} and \eqref{eq:asC}
one directly obtains the upper bound in  \eqref{eq:c1c2}.

We are therefore left with the proof of \eqref{eq:asclaim}.
This follows by observing first that 
$$
Z_{l+r, \theta^{-l}\bgo}^\bv \ge Z_{l, R\bgo}^\bv  Z_{r, \bgo}^\bv e^{\tilde \lambda(\gto_0-\gto_{-l})}, 
$$ 
where $R\bgo$ is the disorder sequence reflected around the 
origin: $(R \bgo)_n= \bgo_{-n}$.
Therefore, with $s=l+r$, we have
\begin{multline}
\label{eq:aslong}
\frac 1 s \log G_{l,r}\left(\bgo \right)
\, \le \, \frac 1s \log K(s)+
\frac{\tilde \gl \left( \tilde \go_r-\gto_0+\gto_{-l} + \tilde h \right)}{s}
+ \\
\frac 1 s \log \left(\frac{ 1+ \exp \left( -2 \gl \sum_{n=-l+1}^r
(\go_n+h)\right)}2\right)
-\frac 1 s \log Z_{l, R\bgo}^\bv - \frac 1 s \log Z_{r, \bgo}^\bv.
\end{multline}
The leading terms are the last two:
by  definition of $\tf(\bv)$, for every $\gep>0$ there exists $n_0(\bgo)$,
$\bbP(\dd\bgo)$--almost surely finite,
such that for $l$ and $r$ larger than $n_0(\bgo)$
both $ (1/ l) \log Z_{l, R\bgo}^\bv $
and $(1 /r)\log Z_{r, \bgo}^\bv$
are larger than $ \tf (\bv) -\gep/2$ and this
easily yields the existence of $s_1(\bgo)$ such that
\begin{equation}
\frac 1 s \log Z_{l, R\bgo}^\bv + \frac 1 s \log Z_{r, \bgo}^\bv
\, \ge  \tf (\bv) -\gep/2,
\end{equation}
for every $s\ge s_1(\bgo)$.
The remaining  term
in the last line of \eqref{eq:aslong} is treated by an analogous
splitting of the sum and by applying the law of large numbers
one sees that it gives a negligible contribution for $s$
sufficiently large. Finally the first  two terms
in the right--hand side of \eqref{eq:aslong}
are both vanishing as $s $ diverges by the properties of $K(\cdot)$,
cf. \eqref{eq:alpha}, and by the fact that $\tilde \go_1$ is integrable,
so that $\tilde \go _n /n \stackrel{n\to \infty}{\longrightarrow} 0$, $\bbP (\dd \bgo)$--a.s..

\hfill $\stackrel{\text{Proposition\; \ref{th:as},\; upper bound}}{\Box}$

\noindent
{\it Lower bound.}
By selecting in the sum in
 \eqref{eq:id1} only the {\sl excursion} from $k$ to 
 $k+s$ we have that we can write for $k+s\le N$
\begin{equation}
  \label{eq:id1min}
  \bP_{N, \bgo}\left( \boldsymbol \Delta_N (k) = s\right) \, \ge
  \,
  K(s)
\frac{Z^{\bv}_{k, \bgo}
\left(
1+ e^{ -2 \gl \sum_{n=k+1}^{k+s} (\go_n+h)}\right)\zeta(\gto_{k+s})Z^{\bv}_{N-k-s,\theta^{k+s}\bgo}}
{2Z^{\bv}_{N,\bgo}},
\end{equation}
and we will use also
\begin{eqnarray}
  \label{eq:id2}
  \frac1{Z^{\bv}_{N,\bgo}}=\frac1{Z_{N,\bgo}(O_{k,s})}(1-\bP_{N,\bgo}(O^{\complement}_{k,s})),
\end{eqnarray}
where $O_{k,s}$ is the event
\begin{eqnarray}
  \label{eq:Olr}
  O_{k,s}=\left\{\exists\; n\in \{k-s^2,\ldots,k+s^2\}:
\; S_{n}=0\right\}
\end{eqnarray}
(with  the conventions of Section \ref{sec:conv} for $Z_{N,\bgo}(\bar\Omega)$).
Thanks to Lemma \ref{th:lbpk}, one can write for some $c>0$
\begin{equation}
  Z_{N,\bgo}(O_{k,s})\le c\, s^c
  Z^{\bv}_{k, \bgo}Z^{\bv}_{s,\theta^{k}\bgo}
Z^{\bv}_{N-k-s,\theta^{k+s}\bgo}\zeta(|\gto_{k}|)\zeta(|\gto_{k+s}|)
\end{equation}
which, together with  \eqref{eq:id1min} and  \eqref{eq:id2},  implies
\begin{eqnarray}
\label{eq:strem}
  \bP_{N, \bgo}\left( \boldsymbol \Delta_N (k) \, =\, s\right)\ge c \, s^{-c}
G_{0,s}(\theta^k\bgo)\left[\zeta(|\gto_{k}|)\zeta(|\gto_{k+s}|)\right]^{-1}
 -\bP_{N,\bgo}(O^{\complement}_{k,s}).
\end{eqnarray}
Proceeding in analogy with the proof of the upper bound, and using 
in addition Lemma \ref{th:lbpk} 
to bound $Z^\bv_{s,\bgo}$ {\it above}, 
one can show that 
$G_{0,s}(\bgo) \exp( (\tf(\bv)+\gep) s )$ diverges $\bbP(\dd \bgo)$--a.s.
as $s\to \infty$,
for every $\gep>0$. 
Therefore, keeping in mind that $k+s \le N$,
the first term on the r.h.s. of \eqref{eq:strem}
is bounded {\it below} by  $\exp(-(\tf(\bv)+\gep)s)$ if $s$ is larger than $s_0(\theta ^k \bgo)$, where 
$s_0(\bgo)$ is 
 a suitable $\bbP(\dd \bgo)$--a.s.
finite number,
and the same quantity may then be bounded below by
\begin{eqnarray}
  C(\theta^k\bgo)\exp(-(\tf(\bv)+\gep)s),
\end{eqnarray}
where  $C(\bgo)$ is a  $\bbP(\dd \bgo)$--a.s. positive random variable.
On the other hand, from the upper bound in Proposition~\ref{th:as} one easily obtains 
\begin{eqnarray}
  \label{eq:s2}
  \bP_{N,\bgo}(O^{\complement}_{k,s})\le \,c C^{(2)}_\gep(\theta^k\bgo)\exp(-(\tf(\bv)-\gep)s^2).
\end{eqnarray}
Taking $\gep$ small enough, this immediately implies the lower bound in Proposition  \ref{th:as}
for $0 \le k \le N-s$.

The restrictions on $k$ is of course due to having chosen in the first step of the proof
\eqref{eq:id1min} the $S$--{\sl excursion} from $k$ to $k+s$.
But we may as well choose the  $S$--{\sl excursion} from $k-s+1$ to $k+1$,
and the argument may be repeated yielding the same bound, except for 
different $\bgo$ dependent constants, for every $k$ ranging from $s-1$ to $N-1$.
Therefore the proof is complete for $k$ in
$\{0, 1, 2, \ldots, N-s\} \cup \{s-1, s, \ldots, N-1\}$ and this is the whole
set of sites smaller than $N$ if $s<N/2$.

\hfill $\stackrel{\text{Proposition\; \ref{th:as},\; lower bound}}{\Box}$

\section{Finite--size corrections to the infinite--volume free energy}

\noindent
{\it Proof of Proposition \ref{th:Nfin1}}.
Just note that
\begin{equation}
\begin{split}
  0\, \le\,  \frac1{2N}\bbE \log Z^{\bv}_{2N,\bgo}- \frac1{N}\bbE \log Z^{\bv}_{N,\bgo}
  \, &=\, 
-\frac1{2N}\bbE \log \bP^{\bv}_{2N,\bgo}(S_N=0)
\\
&\le\, c \frac{\log N}N+\frac{\tilde\gl}N\bbE|\gto_N+\tilde h|,
\end{split}
\end{equation}
where the first inequality is a consequence of the renewal property of $S$ and the last one
of Lemma \ref{th:lbpk}. Inequality \eqref{eq:Nfin1} then immediately follows.

\hfill $\stackrel{\text{Proposition \ref{th:Nfin1}}}{\Box}$

\medskip

\noindent
{\it Proof of Theorem \ref{th:Nfin2}}. It is sufficient to prove that there exists $c(\bv)$ such that, 
for every $N\in\N$,
\begin{eqnarray}
\label{eq:sufficient}
  \frac 1{2N} \bbE\log Z^{\bv}_{2N,\bgo}- \frac 1{N} \bbE\log Z^{\bv}_{N,\bgo}\le \frac{c(\bv)}N.
\end{eqnarray}
It is convenient to define, for $x\ge 0$,
\begin{eqnarray}
  \label{eq:hatF}
  \hat F(\bv,x)=\lim_{N\to\infty}\frac1N \bbE\log \hat Z^{\bv,x}_{N,\bgo}:=
\lim_{N\to\infty}\frac1N \bbE\log \bE \left(e^{\mathcal H_{N,\bgo}(S)+x\sum_{n=1}^N \delta_n}\delta_N\right)
\end{eqnarray}
which, in view of \eqref{eq:normalisaz}, just corresponds to $\hat
{\tf}(\bv,x):=\tf(\gl,h,\tilde \gl, \tilde h+ x/\tilde \gl)$ if
$\tilde\gl\ne0$.  Note that, since $\bv\in\cL$, one has $\partial_{x}
\left. \hat {\tf}(\bv,x)\right|_{x=0}>0$ (cf. Lemma \ref{lemma:derivne0}
below) and that
\begin{equation}
  \label{eq:path}
  \hat {\tf}(\underline 0,x)\ge \tf(\bv)>0,
\end{equation}
provided that  $x\ge(2\gl+\tilde \gl(1+\tilde h))/(\partial_{x}
\left. \hat {\tf}(\bv,x)\right|_{x=0})$. 
Indeed, since $\hat {\tf}$ is convex in $x$
one has 
$$
\hat {\tf}(\bv,x)\ge \hat {\tf}(\bv,0)+x
\partial_{x} \left. \hat {\tf}(\bv,x)\right|_{x=0}
\ge \tf(\bv)+2\gl+\tilde \gl(1+\tilde h)
$$
from which \eqref{eq:path} immediately follows, since $\bbE\, \go_1^2=\bbE\, \gto_1^2=1$.
Essentially the 
same argument shows that there esists a smooth (e.g. differentiable) path $(\bv_t,x_t)$, with $0\le t\le1$, such that
$(\bv_0,x_0)=(\underline 0,x)$, $(\bv_1,x_1)=(\bv,0)$, 
 and such that $\hat {\tf}(\bv_t,x_t)\ge {\tf}(\bv)$ for every $t$.
Note that at $(\underline 0,x)$ the desorder dependence disappears and we have simply a homogenous pinning model 
which, thanks to \eqref{eq:path}, is in the localized phase. For this model it is easy to 
prove (this can be  extracted, for instance, from Appendix A of \cite{cf:GT05}) that 
\begin{eqnarray}
  \label{eq:condiz_iniz}
  \frac 1{2N}\log \hat Z^{\underline 0,x}_{2N,\bgo}- \frac 1{N} \log \hat Z^{\underline 0,x}_{N,\bgo}\le \frac{c_0(\bv)}N.
\end{eqnarray}
On the other hand, we will prove in a moment that 
\begin{eqnarray}
  \label{eq:derivata}
  \left|\frac d{dt} \left( \frac 1{2N} \bbE \log \hat Z^{\bv_t,x_t}_{2N,\bgo}- 
\frac 1{N} \bbE \log \hat Z^{\bv_t,x_t}_{N,\bgo}
\right)\right|\le \frac{c_1(\bv)}N
\end{eqnarray}
uniformly for $0\le t\le 1$. Of course, Eqs. \eqref{eq:condiz_iniz}-\eqref{eq:derivata} immediately give 
\eqref{eq:sufficient}.
To prove \eqref{eq:derivata}, 
let us compute for instance the derivative with respect to $x$: with obvious notations,
\begin{multline}
\left|\frac d{dx} \left( \frac 1{2N} \bbE \log \hat Z^{\bw,x}_{2N,\bgo}- \frac 1{N} \bbE \log \hat Z^{\bw,x}_{N,\bgo}
\right)\right|
\\
=\,  \left|\frac1{2N}\sum_{n=1}^{2N}\bbE \,\hat \bE^{\bw,x}_{2N,\bgo}(\delta_n)-
\frac1{N}\sum_{n=1}^{N}\bbE \,\hat \bE^{\bw,x}_{N,\bgo}(\delta_n)\right|
\end{multline}
which, thanks to  Eq. \eqref{eq:bordo}, is bounded above by
\begin{eqnarray}
  \frac {c_3(\bw,x)}N \sum_{n=1}^N e^{-c_4(\bw,x)(N-n)}\le \frac{c_5(\bw,x)}N,
\end{eqnarray}
$c_5$ being  bounded above uniformly for $(\bw,x)$ belonging to the path. Similar 
estimates hold for the derivatives with respect to $\bw$ and therefore \eqref{eq:derivata} follows.

\hfill $\stackrel{\text{Theorem \ref{th:Nfin2}}}{\Box}$

\section{A Central Limit Theorem}

\label{sec:CLT}
The first step in the proof  of Theorem \ref{th:CLT}
is to show that the variance of $N^{-1/2}\log Z_{N,\bgo}^\bv$ is not trivial in the 
infinite volume limit:
\begin{lemma}
\label{lemma:nontriv}
  If $\bv\in \cL$, then
  \begin{eqnarray}
    \label{eq:var_non_triv}
    0<\liminf_{N\to\infty}\frac1N  {\rm Var} (\log Z_{N,\bgo}^\bv)\le\limsup_{N\to\infty}
\frac1N  {\rm Var} (\log Z^\bv_{N,\bgo})<\infty.
  \end{eqnarray}
\end{lemma}

\medskip

\noindent
{\it Proof of Lemma \ref{lemma:nontriv}.} The upper bound is an immediate consequence of the 
deviation inequality \eqref{eq:Lip} applied to $g(\bgo)=N^{-1/2}\log Z^\bv_{N,\bgo}$.
Indeed, it is immediate to verify that in this case $\Vert g\Vert_{Lip}\le c$, for some constant $c$.

To obtain the lower bound, we employ a martingale method analogous to that developed in \cite{AizWeh90}.
Suppose first that $\tilde\gl\ne0$.
Observe that, if $X_{N,\gto}=\bbE[ \log Z^\bv_{N,\bgo}|\gto]$ and
\begin{equation}
Y^{(N)}_k=\bbE\left[X_{N,\gto}|\gto_1,\ldots,\gto_k\right]-\bbE\left[X_{N,\gto}|\gto_1,\ldots,\gto_{k-1}\right],
\end{equation}
then
\begin{eqnarray}
\label{eq:martingala}
 {\rm  Var} (\log Z^\bv_{N,\bgo})\ge {\rm Var}(X_{N,\gto})=\sum_{k=1}^N \bbE\left[ \left(Y^{(N)}_k\right)^2\right]\ge
\sum_{k=1}^N \bbE\left[\left(\bbE\left(Y^{(N)}_k|\gto_k\right)\right)^2\right],
\end{eqnarray}
with the convention that $\bbE\left[X_{N,\gto}|\gto_1,\ldots,\gto_{k-1}\right]=\bbE\left[ X_{N,\gto}\right]$ if $k=1$.
Next, observe that
\begin{eqnarray}
  \partial_{\gto_k} \bbE\left[Y^{(N)}_k|\gto_k\right]=
\tilde\lambda \,\bbE\left[\bE_{N,\bgo}(\delta_k)|\gto_k\right]\ge0, 
\end{eqnarray}
and that 
\begin{eqnarray}
  \left|\partial^2_{\gto_k} \bbE\left[Y^{(N)}_k|\gto_k\right]\right|\le \tilde\lambda^2.
\end{eqnarray}
At this point, one can employ the following
\begin{lemma} (\cite{AizWeh90})
 Let $0\le M\le1$, $\beta>0$ and let $V_\beta$ be the set of functions 
\begin{equation}
V_\beta=\{g:\mathbb R\to\mathbb R: g\in C^2, 0\le g'(\cdot)\le1, |g''(\cdot)|\le \beta\}
,
\end{equation}
and, for every probability law $\nu$ on $\mathbb R$,
\begin{eqnarray}
  \gamma_\nu(M,\beta)=\inf\left\{\sqrt{\int g^2(\eta) \nu(d\eta)}\left|\,g\in V_\beta, \int g'(\eta) \nu(d\eta)=M\right.
\right\}.
\end{eqnarray}
Then, $\gamma_\nu(M,\beta)$ is convex in $M$ and $\gamma_\nu(M,\beta)>0$ for $M>0$, provided that $\nu$ is not
concentrated on a single point.
\end{lemma}
Identifying $\beta=\tilde\lambda$, $g_k(\gto_k)=\tilde\gl^{-1}\bbE\left[Y^{(N)}_k|\gto_k\right]$ 
and $\nu$ with the law of $\gto_1$, using the convexity of $\gamma_\nu$ and 
recalling \eqref{eq:martingala}, one obtains immediately
\begin{eqnarray}
 {\rm  Var} (\log Z^\bv_{N,\bgo})\ge N{\tilde\gl^2}\left(\gamma_\nu\left(\frac1N\sum_{k=1}^N\bbE\,\bE_{N,\bgo}(
\delta_k),\tilde\lambda\right)\right)^2.
\end{eqnarray}
To conclude the proof of Lemma \ref{lemma:nontriv}, it suffices therefore to show that

\begin{lemma}
\label{lemma:derivne0}
If $\bv\in \cL$,
\begin{eqnarray}
\label{eq:derivne0}
\liminf_{N\to\infty}\frac1N\sum_{k=1}^N\bbE\,\bE_{N,\bgo}(\delta_k)>0.
\end{eqnarray}
\end{lemma}
\medskip

\noindent
{\it Proof of Lemma \ref{lemma:derivne0}.}
It is enough to prove that there exists $\varepsilon$ sufficiently small such that
\begin{eqnarray}
\lim_{N\to\infty}  \bbE\, \bP_{N,\bgo}\left(|\{j: S_j=0\}|\le \varepsilon N\right)=0.
\end{eqnarray}
Indeed, we have
\begin{multline}
\label{eq:stirling}
  \bbE\, \bP_{N,\bgo}\left(|\{j: S_j=0\}|\le \varepsilon N\right)\, =
  \\
  \sum_{\ell=0}^{\lfloor \varepsilon N\rfloor }
\sum_{1\le i_1 <\ldots < i_\ell\le N}  \bbE\, \bP_{N,\bgo}\left(\{j: S_j=0\}=\{i_1,i_2,\ldots,i_\ell\}\right)\\
\le\sum_{\ell=0}^{\lfloor \varepsilon N\rfloor}\left(
\begin{array}{c}
N\\
\ell
  \end{array}
\right)
c_1^{\varepsilon N} e^{-c_2 N}\le c_3(\varepsilon N+1)\left(
\begin{array}{c}
N\\
\lfloor\varepsilon N\rfloor
  \end{array}
\right)
e^{-c_4 N},
\end{multline}
where in the first inequality we employed  \eqref{eq:ritorni}. Using Stirling's approximation, it is clear that,
if $\varepsilon$ is  small enough, the right--hand side of \eqref{eq:stirling} decays exponentially in $N$.

\hfill $\stackrel{\text{Lemma\;  \ref{lemma:derivne0}}}\Box$

In the case $\tilde\gl=0$, the condition $\bv\in \mathcal L$ implies  
$\gl\ne0$ and the proof of Lemma \ref{lemma:nontriv}
proceeds analogously, with the role of $\gto$ being played by
$\go$ and  \eqref{eq:derivne0} replaced by
\begin{eqnarray}
\liminf_{N\to\infty}\frac1N\sum_{k=1}^N\bbE\,\bE_{N,\bgo}(\Delta_k)>0.
\end{eqnarray}
\hfill $\stackrel{\text{Lemma\;  \ref{lemma:nontriv}}}\Box$

The second step in the proof of Theorem \ref{th:CLT} is the observation that
(this follows for instance from Lemma \ref{th:lbpk})
\begin{equation}
\label{eq:taglio}
Z^\bv_{N_1,\bgo}Z^\bv_{N_2,\theta^{N_1}\bgo}\le Z^\bv_{N,\bgo}
\le c N^c e^{\tilde\gl |\gto_{N_1}+\tilde h|}  Z^\bv_{N_1,\bgo}Z^\bv_{N_2,\theta^{N_1}\bgo},
\end{equation}
for some $c<\infty$ independent of $\bgo$.
Therefore, keeping also into account the upper bound in  \eqref{eq:var_non_triv},
one obtains the following approximate subadditivity property:
\begin{equation}
  {\rm Var}(\log Z^\bv_{N,\bgo})\le {\rm Var}(\log Z^\bv_{N_1,\bgo})+{\rm Var}(\log Z^\bv_{N_2,\bgo})+c' \sqrt N \log N,
\end{equation}
for some constant $c'$ and $N$ large enough.
This, together with  \eqref{eq:var_non_triv}, implies that the following limit exists:
\begin{equation}
  \label{eq:varianza}
  \lim_{N\to\infty} \frac1N {\rm Var}(\log Z^\bv_{N,\bgo})=\sigma^2(\bv)>0.
\end{equation}
Next, employing repeatedly   \eqref{eq:taglio}, one obtains the decomposition 
\begin{equation}
  \label{eq:decomposiz}
  \left|\log Z^\bv_{N,\bgo}-\sum_{\ell=0}^{\lfloor N^{1/4}\rfloor-1} \log Z^\bv_{\lfloor N^{3/4}\rfloor,
\theta^{\ell \lfloor N^{3/4}\rfloor}\bgo}\right|\le N^{1/4}\log(c N^c)+
\tilde\gl\sum_{\ell=1}^{\lfloor N^{1/4}\rfloor-1}|\gto_{\ell \lfloor N^{3/4}\rfloor}|
\end{equation}
for every $\bgo$ (we assumed for simplicity that $N$ is multiple of $\lfloor N^{1/4}\rfloor$).
Therefore, if $U_{N,\bgo}=N^{-1/2}(\log Z^\bv_{N,\bgo}- \bbE\, \log Z^\bv_{N,\bgo})$,
\begin{equation}
\label{eq:summand}
U_{N,\bgo}=\frac1{N^{1/8}}\sum_{\ell=0}^{\lfloor N^{1/4}\rfloor -1}
{U_{\lfloor N^{3/4}\rfloor,\theta^{\ell \lfloor N^{3/4}\rfloor}\bgo}}+Q_{N,\bgo},
\end{equation}
where  $Q_{N,\bgo}$
tends to zero in law and $\bbP(\dd \bgo)$--a.s..
The summands in \eqref{eq:summand} are IID random variables and  Lyapunov's condition, which in this case reads simply
$$
\lim_{N\to\infty}\frac1{N^{\delta/8}}\bbE\left[ (U_{\lfloor N^{3/4}\rfloor,\bgo})^{2+\delta}\right]
=0
$$
for some $\delta >0$, follows immediately since the deviation inequality \eqref{eq:Lip} implies that 
$\bbE [(U_{N,\bgo})^{2+\delta}]$ is bounded uniformly in $N$. 
The central limit theorem is proven. 
\hfill $\stackrel{\text{Theorem\; \ref{th:CLT}}}{\Box}$

\appendix

\section{Some technical estimates}
\label{sec:app}

Recall the notations in Section \ref{sec:conv}. 
\medskip

\begin{lemma}
\label{th:lbpk}
There exists $c>0$ 
such that for  every $k,N\in \N$, $k\le N$
 and every $\bgo$
we have that
\begin{equation}
\label{eq:lbpk}
\bP_{N, \bgo}\left( S_k=0\right) \, \ge \, 
\frac c {(k\wedge (N-k))^c}\zeta(|\gto_k|)^{-1}.
\end{equation}
\end{lemma}
\medskip

\noindent
{\it Proof.}
We write
\begin{equation}
\bP_{N, \bgo}\left( S_k=0\right) \,=\,
\frac{
Z^{\bv}_{k, \bgo}Z^{\bv}_{N-k, \theta^k \bgo}}
{Z^{\bv}_{N,\bgo}}
\end{equation}
and 
\begin{multline}
  \label{eq:decompos_Z}
Z^{\bv}_{N,\bgo}=Z^{\bv}_{k, \bgo}Z^{\bv}_{N-k, \theta^k \bgo}\, +
\\
\sumtwo{0\le j<k}{k<\ell\le N}  
Z^{\bv}_{j,\bgo}Z^{\bv}_{N-\ell,\theta^\ell\bgo}\zeta(\gto_\ell)K(\ell-j)\frac{1+e^{-2\gl\sum_{n=j+1}^\ell
(\go_n+h)}}2.
\end{multline}
Using \eqref{eq:alpha} and the definition of slow variation for $L(\cdot)$, one finds that
\begin{eqnarray}
\frac{ K(\ell-j)}{K(k-j)K(\ell-k)}\le c((\ell-k)\wedge (k-j))^{2\alpha}\le c(k\wedge (N-k)))^{2\alpha},
\end{eqnarray}
uniformly in $j,\ell$, 
so that the expression \eqref{eq:decompos_Z} can be bounded above by
\begin{eqnarray}
  c' Z^{\bv}_{k, \bgo}Z^{\bv}_{N-k, \theta^k \bgo}\zeta(|\gto_k|)(k\wedge (N-k)))^{2\alpha},
\end{eqnarray}
for some constant $c'$ independent of $\bgo$, which directly yields \eqref{eq:lbpk}. 

\hfill $\stackrel{\text{Lemma\; \ref{th:lbpk}}}{\Box}$

\begin{lemma}
\label{prop:spezz}
  Let $m\in \N$, $0\le i_0< i_1\ldots < i_m\le N$ with $i_\ell\in \N$, and let $A_j,$ $j=0,\ldots,m-1$ be events
depending only on $S_{i_{j}+1},\ldots,S_{i_{j+1}-1}$, i.e., $A_j\in\sigma(S_n: {i_{j}}< n < {i_{j+1}})$.
Then,
\begin{equation}
  \label{eq:spezz}
  \bP_{N,\bgo}\left(S_{i_j}=0\;\text{\rm for}\; 0\le j\le m; \cap_{j=0}^{m-1} A_j\right)\le
\prod_{j=0}^{m-1}
\bP_{i_{j+1}-i_{j},\theta^{i_{j}}\bgo}(A_j).
\end{equation}
\end{lemma}

\medskip

\noindent
{\it Proof.} 
This is elementary: just rewrite the probability in  \eqref{eq:spezz} as
\begin{equation}
  \label{eq:element}
\frac{\bE\left(e^{\mathcal H_{N,\bgo}(S)}\ind_{\{S_{i_j}=0\;\text{\rm for}\; 0\le j\le m\}}
\ind_{\{\cap_{j=0}^{m-1} A_j\}}\ind_{\{S_N=0\}}\right)}{Z^\bv_{N,\bgo}},
\end{equation}
and notice that one obtains an upper bound for it if one constrains the walk to touch zero at $i_0,\ldots,i_m$ in the
denominator $Z^\bv_{N,\bgo}$. At that point, the probability factorizes thanks to the renewal property of the random
walk and one obtains just the right--hand side of 
 \eqref{eq:spezz}.

\hfill $\stackrel{\text{Lemma\; \ref{prop:spezz}}}{\Box}$
\medskip

\section{Some remarks on $\mu(\bv)$}

\subsection{}
\label{sec:mu}
In this section we sketch a proof of the fact that, under some reasonable 
conditions on the law $\bbP$, the strict inequality $\mu(\bv)< \tf(\bv)$ holds
in the localized region. 

Since $\gl^2+\tilde\gl^2>0$ in $\cL$, let us assume for definiteness
that $\tilde\gl>0$ (in the alternative case, the role of $\gto$ in the following is played by $\go$).
In analogy with \cite{cf:GT05}, we 
assume that one of the  following two conditions holds for the IID sequence  $\gto$:
\newcounter{Lcount}
\begin{list}{{\bf C\arabic{Lcount}}:}
{\usecounter{Lcount}
 \setlength{\rightmargin}{\leftmargin}}
\item {\it Continuous random variables}. 
The law of $\gto_1$
has a  density $P(\cdot)$ with respect to the Lebesgue measure on $\R$, and the function
\begin{eqnarray}
\label{eq:exentropia}
  I: \R\ni x\longrightarrow I(x)=\int_{\R} P(y)\log P(y+x) \dd y ,
\end{eqnarray}
exists and is at least twice differentiable in a neighborhood of  $x=0$. 
This is true in great generality whenever $P(\cdot)$ is positive, for instance in the case of
 $P(\cdot)=\exp(-V(\cdot))$, with $V(\cdot)$  a polynomial bounded below.
\item {\it Bounded random variables}. The random variable $\gto_1$ is bounded,
\begin{eqnarray}
  \label{eq:omegabound}
  |\gto_1|\le M<\infty.  
\end{eqnarray}
\end{list}
Assume first  that condition {\bf C1} holds. Given $\varepsilon>0$, let
 $\tilde\bbP_N$ be the law obtained from
$\bbP$ shifting the law of $\gto_1,\ldots,\gto_N$ so that $\tilde\bbE_N [\gto_\ell]=-\varepsilon$. If $\varepsilon$
is small enough, thanks to the assumed regularity of the function $I(\cdot)$ in \eqref{eq:exentropia} one has
\begin{eqnarray}
  \label{eq:entropia}
  {\rm H}(\tilde \bbP_N|\bbP):=\bbE \left( \frac{\dd\tilde\bbP_N}{\dd\bbP}\log \frac{\dd\tilde\bbP_N}{\dd\bbP}\right)\le 
K N \varepsilon^2,
\end{eqnarray}
for some finite constant $K$.
Then, applying the Jensen inequality, 
\begin{multline}
      \frac1N\log\bbE\left[ \frac {1+e^{-2\lambda\sum_{n=1}^N(\go_n+h)}} {Z_{N, \bgo}^\bv} \right]\,=
      \\
  \frac1N\log\tilde\bbE_N 
\left[ \frac {\left(1+e^{-2\lambda\sum_{n=1}^N(\go_n+h)}\right)} 
{Z_{N, \bgo}^\bv} e^{\log\left({\dd \bbP}/{\dd \tilde \bbP_N}\right)} \right]
\\
\ge \, \frac 1N\tilde\bbE_N \left( \log \frac{\dd \bbP}{\dd \tilde \bbP_N}\right)-\frac1N\bbE\,\log Z^{\bv'}_{N,\bgo},
\end{multline}
so that
\begin{eqnarray}
\label{eq:entrop}
  \mu(\bv)\le K\varepsilon^2+\tf(\bv'),
\end{eqnarray}
where $\bv'$ is obtained from $\bv$ replacing $\tilde h$ with $\tilde h-\varepsilon$.
Since $\tf(\cdot)$ 
is smooth in the localized region and the derivative $\partial_{\tilde h}\tf(\bv)\ge0$ is not zero (cf. 
\eqref{eq:derivne0}), for $\varepsilon$ sufficiently small one has
\begin{eqnarray}
\label{eq:sketch}
  \mu(\bv)\le K\varepsilon^2+\tf(\bv)-
\frac12\partial_{\tilde h}\tf(\bv)\varepsilon\le \tf(\bv)-\frac14\partial_{\tilde h}\tf(\bv)\varepsilon< \tf(\bv).
\end{eqnarray}

In the case of condition {\bf C2}, the proof of  \eqref{eq:sketch} is slightly more complicated. 
Instead of shifting the law of $\gto_1,\ldots,\gto_N$, $\tilde\bbP_N$ is obtained by tilting it:
\begin{eqnarray}
  \label{eq:tildebbP}
 \frac{ \dd \tilde \bbP_N}{\dd \bbP}( \bgo)\,\propto\,
 \exp\left(-\varepsilon \sum_{n=1}^N \gto_n\right).
\end{eqnarray}
The estimate \eqref{eq:entrop} on the relative entropy is still valid, with a different constant $K$, while 
the proof that
\begin{equation}
  \lim_{N\to\infty}\frac1N\tilde\bbE_N \log Z_{N,\bgo}^\bv\le \tf(\bv)-c\varepsilon, 
\end{equation}
for some positive $c$, although rather intuitive in view of the previous case (note in fact that, for $\varepsilon$ small,
the tilting \eqref{eq:tildebbP} produces a shift of  $-\gep+O(\gep^2)$ of the average of $\gto_1,\ldots,\gto_N$), 
still requires some care.
We do not report  here other details, which the interested reader  can  easily reconstruct  from the proof of 
Lemma 3.3 in \cite{cf:GT05}.

\hfill $\Box$

\subsection{The case $\go_1\sim -\go_1$}
\label{sec:mumu}
Here we prove that, if the law of $\go_1$ is symmetric, then the definition \eqref{eq:mu} of $\mu(\bv)$ is equivalent to
\eqref{eq:mu2}.
In fact,
\begin{multline}
  \label{eq:mumu}
  \bbE\frac {e^{-2\gl\sum_{n=1}^N(\go_n+h)}}{Z^\bv_{N,\bgo}}
  \\
  =\, 
\bbE\frac 1{\bE\left[\exp\left(2\gl\sum_{n=1}^N(\go_n+h)(1-\Delta_n)+\tilde\gl\sum_{n=1}^N(\gto_n+\tilde h)\delta_n\right)
\delta_N\right]}\\
\stackrel{S\stackrel{\bP}{\sim}-S}=
\bbE\frac 1{\bE\left[\exp\left(2\gl\sum_{n=1}^N(\go_n+h)\Delta_n+\tilde\gl\sum_{n=1}^N(\gto_n+\tilde h)\delta_n\right)
\delta_N\right]}\\
\le
\bbE
\frac 1{\bE\left[\exp\left(2\gl\sum_{n=1}^N(\go_n-h)\Delta_n+\tilde\gl\sum_{n=1}^N(\gto_n+\tilde h)\delta_n\right)
\delta_N\right]}\stackrel{\omega\stackrel{\bbP}{\sim}-\go}=\bbE\frac 1{Z^\bv_{N,\bgo}},
\end{multline}
where in the inequality we simply used the fact that $\gl,h\ge0$, see Section \ref{sec:model}.
Therefore, 
\begin{equation}
\begin{split}
  \frac 1N \log \bbE
\left[ \frac 1 {Z_{N, \bgo}^\bv} \right]\, &\le\, 
\frac 1N \log \bbE
\left[ \frac {1+e^{-2\lambda\sum_{n=1}^N(\go_n+h)}} {Z_{N, \bgo}^\bv} \right]
\\
&\le \, \frac 1N \log \bbE
\left[ \frac 1 {Z_{N, \bgo}^\bv} \right]+\frac{\log 2}N,
\end{split}
\end{equation}
and the claim follows in the limit $N\to\infty.$

\section*{Acknowledgments}  We are very grateful to an anonymous 
referee for his observations  and, in particular, for having pointed out that our decay of correlation 
estimates yield   
the almost sure result of Th.~\ref{th:correlazioni} (formula \eqref{eq:exlim}).
We would like to thank also 
Francesco Caravenna 
for interesting discussions on the content of Section 5.
This research has been conducted in the framework of the GIP--ANR project  JC05\_42461
({\sl POLINTBIO}).

\bibliographystyle{alea2}

\bibliography{}

\begin{thebibliography}{15}  

\bibitem[Aizenman and Wehr(1990)]{AizWeh90}
M. Aizenman, J. Wehr, \textit{Rounding effects of quenched randomness on first-order phase transitions}, 
Comm. Math. Phys.
 {\bf 130} (1990),  489--528. 

\bibitem[Albeverio and Zhou(1996)]{cf:AZ}
S. Albeverio and X. Y. Zhou, \textit{Free energy and some sample path properties of a random walk with random potential},
  J. Statist. Phys. {\bf 83 } (1996),  573--622.
  

 \bibitem[Alexander and Sidoravicius(2005)]{cf:AS}
 K. S. Alexander and V. Sidoravicius,
 \textit{Pinning of polymers and interfaces by random potentials},
  arXiv e-Print archive: math.PR/0501028.

  
 \bibitem[Biskup and den Hollander(1999)]{cf:BisdH} 
M.  Biskup and F. den Hollander, \textit{A heteropolymer near a linear interface},  Ann. Appl. Probab. {\bf 9} (1999), 668--687.


\bibitem[Bodineau and Giacomin(2004)]{cf:BG}
T. Bodineau and G. Giacomin, \textit{On the localization transition   of random copolymers near selective interfaces},
J. Statist. Phys. {\bf 117} (2004), 801-818.

\bibitem[Bolthausen and den Hollander(1997)]{cf:BdH} 
E. Bolthausen and F. den Hollander,    
 \textit{Localization transition for a polymer near an interface},   
 Ann. Probab.  {\bf 25}  (1997),  1334--1366.   
   
   \bibitem[Caravenna et al.(2005)]{cf:CGG} 
   F. Caravenna, G. Giacomin and M. Gubinelli,
\textit{A numerical approach to copolymers at selective interfaces},
arXiv e-Print archive:  math-ph/0509065, accepted for
publication on J. Statist. Phys.. 


\bibitem[Derrida et al.(1992)]{cf:DHV}
B. Derrida, V. Hakim and J. Vannimenus,
\textit{Effect of disorder on two--dimensional wetting},
J. Statist. Phys. {\bf 66} (1992), 1189--1213.

\bibitem[von Dreifus et al.(1995)]{cf:griffiths} 
H. von Dreifus, A. Klein, J. Fernando Perez,  \textit{Taming 
Griffiths' singularities: infinite differentiability of quenched correlation functions},
  Comm. Math. Phys.  {\bf 170} (1995), 21--39.

\bibitem[Feller(1971)]{cf:Feller2}  
W.~Feller, \textit{An introduction to probability theory and its applications}, Vol. II,  
$2^{\textrm{nd}}$ edition, John Wiley \& Sons, Inc.,   
New York--London--Sydney, 1971.  

\bibitem[Forgacs et al.(1986)]{cf:FLNO}
G. Forgacs, J. M. Luck, Th. M. Nieuwenhuizen and H. Orland,
\textit{Wetting of a disordered substrate: exact critical behavior in two dimensions},
Phys. Rev. Lett. {\bf 57} (1986), 2184--2187.


\bibitem[Garel et al.(1989)]{cf:GHLO}
T. Garel, D. A. Huse, S. Leibler and H. Orland,
\textit{Localization transition of random chains at interfaces},
Europhys. Lett. {\bf 8} (1989), 9--13.

\bibitem[Giacomin(2004)]{cf:G}  
G.~Giacomin, \textit{Localization phenomena in random polymer models}, preprint (2004). 
Available
online: \\
http://www.proba.jussieu.fr/pageperso/giacomin/pub/publicat.html

\bibitem[Giacomin and Toninelli(2005a)]{cf:GT} 
G.~Giacomin, F.~L.~Toninelli, \textit{Estimates on path delocalization for copolymers 
at selective interfaces}, Probab. Theor. Rel. Fields {\bf 133}  (2005), 464--482.

\bibitem[Giacomin and Toninelli(2005b)]{cf:GT05} G.~Giacomin, F.~L.~Toninelli, \textit{Smoothing effect of quenched disorder on 
polymer depinning transitions}, arXiv e-Print archive: math.PR/0506431, accepted for
publication on
Commun. Math. Phys..

\bibitem[Giacomin and Toninelli(2006)]{cf:GTlett} G.~Giacomin, F.~L.~Toninelli, \textit{Smoothing of 
Depinning Transitions for Directed Polymers with Quenched Disorder}, 
Phys. Rev. Lett. {\bf 96} (2006), 070602.

\bibitem[Ledoux(2001)]{cf:LedouxAMS}
M.~Ledoux,
\textit{The concentration of measure phenomenon}, Mathematical Surveys and Monographs, {\bf 89}, American Mathematical Society, Providence, RI, 2001.

\bibitem[Ledoux(2003)]{cf:Ledoux}
M.~Ledoux,
\textit{
Measure concentration, transportation cost, and functional inequalities},
Summer School on Singular Phenomena and Scaling in Mathematical Models, Bonn, 10--13 June 2003.
Available online: http://www.lsp.ups-tlse.fr/Ledoux/  

\bibitem[Monthus(2000)]{cf:Monthus} C. Monthus, \textit{On the
    localization of random heteropolymers at the interface between two
    selective solvents}, Eur. Phys. J. B {\bf 13} (2000), 111--130.
  
\bibitem[Petrelis(2005)]{cf:Petrelis}
N. Petrelis,
\textit{Polymer pinning at an interface},
 arXiv  e-Print archive: math.PR/ 0504464.


\bibitem[Sinai(1993)]{cf:Sinai}
Ya. G. Sinai, \textit{A 
random walk with a random potential},  Theory Probab. Appl. {\bf 38} (1993),  382--385.

\bibitem[Soteros and Whittington(2004)]{cf:SW}
C. E. Soteros and S. G. Whittington, \textit{The statistical mechanics of random copolymers},
 J. Phys. A: Math. Gen. {\bf 37} (2004), R279--R325.

\bibitem[Talagrand(1996)]{cf:newlook} 
M. Talagrand, \textit{A new look at independence}, Ann. Probab. {\bf 24}  (1996), 1--34.

\bibitem[Villani(2003)]{cf:cedric} C. Villani,
\textit{Topics in optimal transportation}, Graduate Studies in Mathematics {\bf 58}, 
American Mathematical Society, Providence, RI, 2003.





\end{thebibliography}

\end{document}